\documentclass[12pt]{article}
\usepackage[margin=1in]{geometry}

\usepackage{mathtools,amsthm,mathrsfs,amssymb}

\newtheorem{dfn}{Definition}[section]
\newtheorem{theorem}[dfn]{Theorem}
\newtheorem{lemma}[dfn]{Lemma}

\newtheorem{corollary}[dfn]{Corollary}
\newtheorem{question}[dfn]{Question}
\newtheorem{conjecture}[dfn]{Conjecture}
\newenvironment{pf}{\noindent{\bf Proof.}}
{\enspace\vrule height5pt depth0pt width5pt} \def\deg{{\rm deg}}

\newtheorem{clm}{Claim}

\makeatletter
\newenvironment{claimproof}[1][\proofname]{\par\pushQED{\hfill$\lozenge$}\normalfont \topsep6\p@\@plus6\p@\relax \trivlist \item\relax{\itshape#1\@addpunct{.}}\hspace\labelsep\ignorespaces}{\popQED\endtrivlist\@endpefalse}
\makeatother

\def\F {{\mathcal F}}

\def\P {{\mathcal P}}
\def\X {{\mathcal X}}

\def\Y {{\mathcal Y}}

\def\dist {{\rm dist}}

\title{Weak diameter choosability of graphs with an excluded minor}
\author{Joshua Crouch \\ 
\small Department of Mathematics, \\
\small Texas A\&M University,\\
\small College Station, TX 77843-3368, USA \\
\and Chun-Hung Liu\thanks{chliu@tamu.edu. Partially supported by NSF under award DMS-1954054 and CAREER award DMS-2144042.} \\
\small Department of Mathematics, \\
\small Texas A\&M University,\\
\small College Station, TX 77843-3368, USA \\
}

\begin{document}

\maketitle

\begin{abstract}
Weak diameter coloring of graphs recently attracted attention, partially due to its connection to asymptotic dimension of metric spaces.
We consider weak diameter list-coloring of graphs in this paper.
Dvo\v{r}\'{a}k and Norin proved that graphs with bounded Euler genus are 3-choosable with bounded weak diameter.
In this paper, we extend their result by showing that for every graph $H$, $H$-minor free graphs are 3-choosable with bounded weak diameter.
The upper bound 3 is optimal and it strengthens an earlier result for non-list-coloring $H$-minor free graphs with bounded weak diameter.
As a corollary, $H$-minor free graphs with bounded maximum degree are 3-choosable with bounded clustering, strengthening an earlier result for non-list-coloring. 

When $H$ is planar, we prove a much stronger result: for every 2-list-assignment $L$ of an $H$-minor free graph, every precoloring with bounded weak diameter can be extended to an $L$-coloring with bounded weak diameter.
It is a common generalization of earlier results for non-list-coloring with bounded weak diameter and for list-coloring with bounded clustering without allowing precolorings.
As a corollary, for any planar graph $H$ and $H$-minor free graph $G$, there are exponentially many list-colorings of $G$ with bounded weak diameter (and with bounded clustering if $G$ also has bounded maximum degree); and every graph with bounded layered tree-width and bounded maximum degree has exponentially many 3-colorings with bounded clustering.

We also show that the aforementioned results for list-coloring cannot be extended to odd minor free graphs by showing that some bipartite graphs with maximum degree $\Delta$ are $k$-choosable with bounded weak diameter only when $k=\Omega(\log\Delta/\log\log\Delta)$.
On the other hand, we show that odd $H$-minor graphs are 3-colorable with bounded weak diameter, implying an earlier result about clustered coloring of odd $H$-minor free graphs with bounded maximum degree.
\end{abstract}

\section{Introduction}\label{sec:intro}
Weak diameter coloring has attracted attention recently in the graph theory community partially due to its connection to asymptotic dimension of metric spaces and clustered coloring of graphs (see \cite{bbeglps}, for example).
We consider weak diameter list-coloring in this paper.

Let $G$ be a graph.
The {\it distance in $G$} between two vertices $a$ and $b$ of $G$, denoted by $\dist_G(a,b)$, is the infimum of the length of a path in $G$ between $a$ and $b$.
For a subset $S$ of $V(G)$, the {\it weak diameter in $G$ of $S$} is $\sup_{a,b \in S}\dist_G(a,b)$. 
For a function $f$ whose domain is a subset of $V(G)$, an {\it $f$-monochromatic component} is a component of the subgraph of $G$ induced by $\{v \in V(G): f(v)=i\}$ for some element $i$ in the image of $f$; for every real number $N$, we say that $f$ has {\it weak diameter in $G$} at most $N$ if $V(C)$ has weak diameter in $G$ at most $N$ for every $f$-monochromatic component $C$.

For a graph $G$ and a positive integer $m$, an {\it $m$-coloring} of $G$ is a function $f: V(G) \rightarrow [m]$.
Note that an $m$-coloring of $G$ is equivalent to a partition of $V(G)$ into at most $m$ parts.
Hence an $m$-coloring of $G$ with small weak diameter in $G$ is a partition of $V(G)$ into at most $m$ sets with small weak diameter.
Such a partition is useful in metric geometry and distributive computing; for example, see the discussion in \cite{bbeglps}.
Moreover, $m$-colorings with weak diameter in $G$ zero are exactly proper $m$-colorings of $G$.

For positive integers $m$ and $N$, we say that a graph $G$ is {\it $m$-colorable with weak diameter in $G$ at most $N$} if there exists an $m$-coloring $f$ of $G$ such that $f$ has weak diameter in $G$ at most $N$.
There are graph classes $\F$ such that all graphs $G$ in $\F$ have $m$-colorings with weak diameter in $G$ at most $N$ for some absolute constant $m$ and for some integer $N$ only dependent on $\F$ \cite{bbeglps,l_asdim}, even though there is no absolute constant $m'$ such that all graphs in $\F$ have proper $m'$-colorings or even just $m'$-colorings whose every monochromatic component has bounded diameter (where the distance is computed in the monochromatic component instead of in the whole graph) (see \cite{lo}).

We mainly consider list-coloring of graphs with small weak diameter in minor-closed families in this paper.

Let $G$ be a graph.
A {\it list-assignment} of $G$ is a function $L$ with domain $V(G)$ such that $L(v)$ is a nonempty set for every $v \in V(G)$.
An {\it $L$-coloring} of $G$ is a function $f$ with domain $V(G)$ such that $f(v) \in L(v)$ for every $v \in V(G)$.
Let $m$ be a positive integer.
An {\it $m$-list-assignment} of $G$ is a function $L$ with domain $V(G)$ such that $|L(v)|=m$ for every $v \in V(G)$.
For a real number $N$, we say that $G$ is {\it $m$-choosable with weak diameter in $G$ at most $N$} if for every $m$-list-assignment $L$ of $G$, there exists an $L$-coloring of $G$ with weak diameter in $G$ at most $N$.

It is known that graphs embeddable in a fixed surface can be 3-colored with small weak diameter \cite{bbeglps,bbegps}.
Dvo\v{r}\'{a}k and Norin \cite{dn} proved that the same result holds even for list-coloring.

\begin{theorem}[{{\cite[Theorem 5(i)]{dn}}}] \label{surface_ch}
For every surface\footnote{In this paper, a {\it surface} is a connected 2-dimensional compact manifold with no boundary.} $\Sigma$, there exists an integer $N$ such that every graph that can be drawn in $\Sigma$ is 3-choosable with weak diameter in $G$ at most $N$.
\end{theorem}

The number of colors in Theorem \ref{surface_ch} is optimal even when $\Sigma$ is the sphere by the Hex Lemma \cite{g} stating that every 2-coloring of a sufficiently large triangular grid has a long monochromatic path.

A graph $H$ is a {\it minor} of another graph $G$ if $H$ is isomorphic to a graph that can be obtained from a subgraph of $G$ by contracting edges.
We say that $G$ is {\it $H$-minor free} if $H$ is not a minor of $G$.

Graphs that can be drawn in a fixed surface are typical examples of graphs with an excluded minor.
For every surface $\Sigma$, there exists an integer $t$ such that the complete graph $K_t$ cannot be drawn in $\Sigma$ by Euler's formula; and every graph that is obtained by contracting edges from a subgraph of a graph drawn in $\Sigma$ can still be drawn in $\Sigma$; so every graph that can be drawn in $\Sigma$ is $K_t$-minor free.

Graphs with an excluded minor are much more general than graphs drawn in surfaces.
In fact, for every graph property $P$ that is not satisfied by all graphs but is preserved under vertex-deletion, edge-deletion and edge-contraction, there exists a graph $H$ such that every graph satisfying $P$ is $H$-minor free.
To name a few such properties $P$: linkless embeddability, knotless embeddability, having a vertex-cover or a feedback vertex set with size at most a fixed constant, and having Colin de Verdi\`{e}re parameter at most a fixed constant \cite{c}.

\subsection{Main results}

\subsubsection{With an excluded minor}

The first main result in this paper is the strengthening of Theorem \ref{surface_ch} to all graphs with an excluded minor.

\begin{theorem} \label{minor_intro}
For every graph $H$, there exists an integer $N$ such that every $H$-minor free graph is 3-choosable with weak diameter in $G$ at most $N$.
\end{theorem}

Recall that the Hex Lemma states that triangular grids require three colors for any coloring of bounded weak diameter.
So the bound 3 is optimal in Theorem \ref{minor_intro} whenever $H$-minor free graphs contains all triangular grids; in particular, it is the case when $H$ is non-planar.
We will show in Theorem \ref{tw_intro} below, the number 3 can be improved to 2 when $H$ is planar.

Theorem \ref{minor_intro} can be applied to obtain upper bounds for the number of colors for clustered coloring as follows.

Let $G$ be a graph.
Let $N$ be an integer.
For every function $f$ with domain $V(G)$, we say that $f$ has {\it clustering} at most $N$ if every $f$-monochromatic component contains at most $N$ vertices.
We say that $G$ is {\it $m$-choosable with clustering at most $N$} if for every $m$-list-assignment $L$ of $G$, there exists an $L$-coloring of $G$ with clustering at most $N$.

Clustered coloring has been extensively studied.
For example, the second author and Wood \cite{lw_hajos} proved that $K_{t+1}$-topological minor free graphs are $(4t-5)$-colorable with bounded clustering, showing that a linear upper bound exists for clustered coloring in constrast to the $\Omega(t^2/\log t)$ lower bound for proper coloring \cite{ef}; Dujmovi\'{c}, Esperet, Morin and Wood \cite{demw} proved that $K_{t+1}$-minor free graphs are $t$-colorable with bounded clustering, showing that the clustered coloring relaxation of Hadwiger's conjecture holds; in fact, Dvo\v{r}\'{a}k and Norin \cite{dn_island} announced in an earlier paper that a forthcoming paper will prove that $K_{t+1}$-minor free graphs are $t$-choosable with bounded clustering.
See \cite{w} for a survey about clustered coloring.

Every coloring of a graph $G$ with clustering at most an integer $N$ has weak diameter in $G$ at most $N$.
The converse asymptotically holds when $G$ has bounded maximum degree, since if $G$ has maximum degree at most an integer $\Delta$, then every subset of $V(G)$ with weak diameter in $G$ at most $k$ contains at most $\sum_{i=0}^k\Delta^i \leq (k+1)\Delta^{k}$ vertices.

Therefore, Theorem \ref{minor_intro} immediately implies the following corollary that generalizes a result of the second author and Oum \cite{lo} stating that for any graph $H$ and integer $\Delta$, $H$-minor free graphs with maximum degree at most $\Delta$ are 3-colorable with bounded clustering.

\begin{corollary} \label{cor_minor_clu_intro}
For any graph $H$ and integer $\Delta$, there exists an integer $N$ such that every $H$-minor free graph with maximum degree at most $\Delta$ is 3-choosable with clustering at most $N$.
\end{corollary}

Again, the number 3 in Corollary \ref{cor_minor_clu_intro} is optimal due to the Hex Lemma since triangular grids have maximum degree at most six.
When $H$ is planar, the number 3 in Corollary \ref{cor_minor_clu_intro} can be reduced to 2 as shown in \cite{lw_minor}.
In this paper, we prove a much stronger result that allows precolorings (see Section \ref{subsubsec:planar}).

\subsubsection{With an excluded planar minor} \label{subsubsec:planar}

For a function $f$ and a subset $S$ of the domain of $f$, we denote by $f|_S$ the function obtained by restricting $f$ on $S$; that is, $f|_S$ is the function with domain $S$ with $f|_S(x)=f(x)$ for every $x \in S$.
Let $f_1$ and $f_2$ be functions such that the domain $S_1$ of $f_1$ is a subset of the domain of $f_2$.
We say that $f_1$ {\it can be extended} to $f_2$ if $f_2|_{S_1}=f_1$.

Let $G$ be a graph.
Let $L$ be a list-assignment of $G$.
A {\it pre-$L$-coloring} of $G$ is an $L|_S$-coloring of $G[S]$ for some subset $S$ of $V(G)$.

The following is the second main result of this paper, roughly stating that every precoloring of an $H$-minor free graph with bounded weak diameter can be extended to a coloring with bounded weak diameter, as long as $H$ is planar and two colors are available for each vertex, generalizing a result of \cite{bbeglps,l_asdim} that only works for 2-(non-list-)coloring and essentially does not allow precolored vertices.

\begin{theorem} \label{tw_intro}
For any planar graph $H$ and integer $N_0$, there exists an integer $N$ such that for every $H$-minor free graph $G$, every $2$-list-assignment $L$ of $G$, and every pre-$L$-coloring $c_0$ of $G$ with weak diameter in $G$ at most $N_0$, $c_0$ can be extended to an $L$-coloring of $G$ with weak diameter in $G$ at most $N$.
\end{theorem}

Note that the condition for precolorings in Theorem \ref{tw_intro} is the weakest that one can expect, because only precolorings with bounded weak diameter can be extended to a coloring with bounded weak diameter.
By applying this result to graphs with bounded maximum degree, we obtain the following result for clustered coloring.

\begin{corollary} \label{cor_tw_clu_intro}
For any planar graph $H$ and integers $\Delta,N_0$, there exists an integer $N$ such that for every $H$-minor free graph $G$ with maximum degree at most $\Delta$, every $2$-list-assignment $L$ of $G$, and every pre-$L$-coloring $c_0$ of $G$ with weak diameter in $G$ at most $N_0$, $c_0$ can be extended to an $L$-coloring of $G$ with clustering at most $N$.
\end{corollary}

Corollary \ref{cor_tw_clu_intro} implies a result in \cite{lw_minor} that essentially does not allow precolorings.
Note that the results in \cite{bbeglps,l_asdim,lw_minor} do allow precolorings; however, the precolored vertices in those results are required to either be close to a set of bounded number of vertices or form a set of bounded size; it is easy to show that the existence of a coloring with bounded weak diameter or with bounded clustering is equivalent to the existence of a desired coloring that can be extended from a precoloring with those requirements; so allowing precolored vertices under those requirements does not increase the strength of the results, in contrast to Theorem \ref{tw_intro} and Corollary \ref{cor_tw_clu_intro} that allow precolorings with the weakest possible assumption.

The strength that allows precolorings in Theorem \ref{tw_intro} and Corollary \ref{cor_tw_clu_intro} is of independent interest.
For example, we will use Theorem \ref{tw_intro} and Corollary \ref{cor_tw_clu_intro} to show that there are indeed exponentially many colorings with bounded weak diameter or bounded clustering (Corollary \ref{cor_exponentially_color_intro}).
Before stating it, we first remark a relation between tree-width and Theorem \ref{tw_intro} and Corollary \ref{cor_tw_clu_intro}.

A {\it tree-decomposition} of a graph $G$ is a pair $(T,\X)$, where $T$ is a tree and $\X$ is a collection $(X_t: t \in V(T))$ of subsets of $V(G)$ such that
    \begin{itemize}
        \item $\bigcup_{t \in V(T)}X_t=V(G)$,
        \item for every $uv \in E(G)$, there exists $t \in V(T)$ such that $\{u,v\} \subseteq X_t$, and
        \item for every $v \in V(G)$, the set $\{t \in V(T): v \in X_t\}$ induces a connected subgraph of $T$.
    \end{itemize}
For every $t \in V(T)$, the set $X_t$ is called the {\it bag} at $t$.
The {\it adhesion} of $(T,\X)$ is $\max_{tt' \in E(T)}|X_t \cap X_{t'}|$.
The {\it width} of $(T,\X)$ is $\max_{t \in V(T)}|X_t|-1$.
The {\it tree-width} of $G$ is the minimum width of a tree-decomposition of $G$.

It is easy to show that tree-width is a minor-monotone graph parameter and planar graphs can have arbitrarily large tree-width.
So for every integer $w$, there exists a planar graph $H$ such that every graph with tree-width at most $w$ is $H$-minor free.
Hence Theorem \ref{tw_intro} and Corollary \ref{cor_tw_clu_intro} can be applied to graphs with bounded tree-width.

Moreover, it is well-known that every graph $G$ with tree-width at most $w$ is $w$-degenerate and hence properly $(w+1)$-colorable, so it contains a stable set $S$ with size at least $|V(G)|/(w+1)$.
Note that for every $m$-list-assignment $L$ of $G$, there are $m^{|S|} \geq m^{|V(G)|/(w+1)}$ different pre-$L$-colorings on $S$, and every pre-$L$-coloring on $S$  has weak diameter in $G$ at most 0.
Therefore, Theorem \ref{tw_intro} and Corollary \ref{cor_tw_clu_intro} immediately imply the following corollary about the existence of exponentially many colorings.

\begin{corollary} \label{cor_exponentially_color_intro}
Let $w$ be a nonnegative integer.
Then 
    \begin{enumerate}
        \item there exists an integer $N$ such that for every integer $m \geq 2$ and every $m$-list-assignment $L$ of a graph $G$ with tree-width at most $w$, there exist at least $m^{|V(G)|/(w+1)}$ different $L$-colorings of $G$ with weak diameter in $G$ at most $N$;
        \item for every integer $\Delta$, there exists an integer $N$ such that for every integer $m \geq 2$ and every $m$-list-assignment $L$ of a graph $G$ with tree-width at most $w$, there exist at least $m^{|V(G)|/(w+1)}$ different $L$-colorings of $G$ with clustering at most $N$.
    \end{enumerate}
\end{corollary}

\subsubsection{Layered tree-width}

Corollary \ref{cor_exponentially_color_intro} leads to further applications to graphs with bounded layered tree-width and bounded maximum degree to show the existence of exponentially many clustered 3-colorings, which applies to graphs with certain geometric properties.

A {\it layering} of a graph $G$ is an ordered partition $(V_1,V_2,...)$ of $V(G)$ into possibly empty sets such that for every edge $e$ of $G$, there exists $i \in {\mathbb N}$ such that the ends of $e$ are contained in $V_i \cup V_{i+1}$.
The {\it layered tree-width} of a graph $G$ is the minimum integer $w$ such that there exists a tree-decomposition $(T,(X_t: t \in V(T)))$ of $G$ and a layering $(V_1,V_2,...)$ of $G$ such that $|X_t \cap V_i| \leq w$ for every $t \in V(T)$ and $i \in {\mathbb N}$.

Graphs with bounded layered tree-width include graphs with certain geometric properties and are incomparable with graphs with forbidden minors.
For example, for nonnegative integers $g$ and $k$, every $(g,k)$-planar graph has layered tree-width at most $(4g+6)(k+1)$ \cite{dmw}, where a graph is {\it $(g,k)$-planar} if it can be drawn in a surface of Euler genus at most $g$ such that every edge contains at most $k$ crossings.
In particular, $(g,0)$-planar graphs are exactly the graphs that can be embedded in a surface of Euler genus at most $g$.
In addition, $(0,1)$-planar graphs are also known as 1-planar graphs, which have been extensively studied in the literature.
It is well-known that every graph is a minor of some $(0,1)$-planar graph.
Another example of graphs with bounded layered tree-width is {\it $(g,k)$-string graphs}, which are the intersection graphs of curves in a surface of Euler genus at most $g$ such that every curve contains at most $k$ intersecting points with other curves.
See \cite{dmw} for other examples for graphs with bounded layered tree-width.

The second author and Wood \cite{lw} proved that every graph with bounded layered tree-width and bounded maximum degree is 3-colorable with bounded clustering.
This result was improved in three different directions: the bound for the clustering was improved in \cite{demww}; the condition for having bounded maximum degree was relaxed to forbidden subgraphs \cite{lw_layer_subgraph}; the weak coloring version was proved in \cite{bbeglps,l_asdim}.

By combining with the proof in \cite{lw}, Corollary \ref{cor_exponentially_color_intro} gives another strengthening of the result in \cite{lw} by showing that there are actually exponentially many 3-colorings with bounded clustering.

\begin{corollary} \label{cor_exponentially_color_layer_intro}
For any positive integers $w$ and $\Delta$, there exists a positive integer $N$ such that every graph with layered tree-width at most $w$ and maximum degree at most $\Delta$ has at least $2^{|V(G)|/(3w)}$ different colorings with clustering at most $N$.
\end{corollary}

\begin{pf}
Let $G$ be a graph with layered tree-width at most $w$ and with maximum degree at most $\Delta$.
Let $(V_1,V_2,...)$ be a layering of $G$ and $(T,\X)$ be a tree-decomposition such that the intersection of any layer and any bag has size at most $w$.
For $i \in [3]$, let $U_i = \bigcup_{k \in {\mathbb N} \cup \{0\}}V_{3k+i}$.
Since $U_1 \cup U_2 \cup U_3 = V(G)$, by possibly shifting the indices of the layering, we may assume $|U_1| \geq |V(G)|/3$.
Since $G[U_1]$ has tree-width at most $w-1$, by Statement 2 of Corollary \ref{cor_exponentially_color_intro}, there exist a positive integer $N_1$ (only dependent on $w$ and $\Delta$) and at least $2^{|U_1|/w} \geq 2^{|V(G)|/(3w)}$ different 2-colorings of $G[U_1]$ with clustering at most $N_1$.
The proof of \cite[Theorem 8]{lw} shows that for every integer $N_2$, every 2-coloring of $G[U_1]$ with clustering at most $N_2$ can be extended to a 3-coloring of $G$ with clustering at most $N_3$ for some integer $N_3$ only dependent on $N_2,w,\Delta$.
Hence there are at least $2^{|V(G)|/(3w)}$ different 3-colorings of $G$ with clustering at most $N$ for some integer $N$ only dependent on $w$ and $\Delta$.
\end{pf}

\bigskip

Note that Corollary \ref{cor_exponentially_color_layer_intro} implies that every $(g,k)$-planar graph (in particular, graphs of bounded Euler genus) with bounded maximum degree has exponentially many 3-colorings with bounded clustering.
The list-coloring version of Corollary \ref{cor_exponentially_color_layer_intro} remains open, even just for the existence of such a coloring.
See Conjecture \ref{conj_choosable_layer_deg}.

\subsubsection{Bipartite graphs and odd minors}

Finally, we consider odd minors.
A graph $H$ is an {\it odd-minor} of another graph $G$ if there exists a set $\{X_h: h \in V(H)\}$ of disjoint subgraphs of $G$ isomorphic to trees such that each tree $X_h$ has a proper 2-coloring $c_h: V(X_h) \rightarrow [2]$ such that for every edge $uv$ of $H$, there exists an edge $u'v'$ of $G$ such that $u' \in V(X_u)$, $v' \in V(X_v)$ and $c_u(u')=c_v(v')$.
It is clear that every odd minor of $G$ is also a minor of $G$.
But the converse is not true: it is easy to show that if $G$ is bipartite, then every odd minor of $G$ is bipartite.

We show that Theorem \ref{minor_intro} and Corollary \ref{cor_minor_clu_intro} cannot be strengthened to odd $H$-minor free graphs for non-bipartite $H$ by showing the following result.

\begin{theorem} \label{bipartite_intro}
There are infinitely many positive integers $\Delta$ such that for every positive integer $N$, there exist a bipartite graph $G$ with maximum degree at most $\Delta$ and a $\lfloor \frac{\log\Delta}{12\log\log\Delta} \rfloor$-list-assignment $L$ such that every $L$-coloring of $G$ has weak diameter in $G$ greater than $N$.
\end{theorem}

Theorem \ref{bipartite_intro} implies that if all bipartite graphs with maximum degree at most $\Delta$ are $k$-choosable with bounded weak diameter, then $k=\Omega(\frac{\log\Delta}{\log\log\Delta})$.
So for every non-bipartite graph $H$, there is no integer $k$ (depending on $H$ or not) such that every odd $H$-minor free graph is $k$-choosable with bounded weak diameter. 

Note that Alon and Krivelevich \cite{ak} conjectured that bipartite graphs with maximum degree at most $\Delta$ are $O(\log\Delta)$-choosable (with weak diameter 0).
If this conjecture is true, then the $\Omega(\frac{\log\Delta}{\log\log\Delta})$ lower bound in Theorem \ref{bipartite_intro} is optimal up to a $\log\log\Delta$ factor. 

On the other hand, we show that the non-list-coloring version of Theorem \ref{minor_intro} can be generalized to odd minor free graphs.

\begin{theorem} \label{odd_minor_intro}
For every graph $H$, there exists an integer $N$ such that every odd $H$-minor free graph $G$ has a 3-coloring with weak diameter in $G$ at most $N$.
\end{theorem}

Theorem \ref{odd_minor_intro} immediately implies the following result in \cite{lo}. 

\begin{corollary}[{\cite[Theorem 1.6]{lo}}] \label{cor_odd_minor_clu_intro}
For any graph $H$ and positive integer $\Delta$, there exists an integer $N$ such that every odd $H$-minor free graph $G$ with maximum degree at most $\Delta$ has a 3-coloring with clustering at most $N$.
\end{corollary}

The bound 3 in Theorem \ref{odd_minor_intro} and Corollary \ref{cor_odd_minor_clu_intro} is optimal for non-planar $H$ by the Hex Lemma.

\subsection{Proof sketch and organization of the paper}

The proofs of Theorems \ref{minor_intro}, \ref{tw_intro} and \ref{odd_minor_intro} follow from the same strategy: we shall apply structure theorems for (odd) $H$-minor free graphs to obtain a tree-decomposition such that every ``torso'', which is a graph obtained from the subgraph induced by a bag by adding extra vertices and edges in a certain way, has certain nice properties, then show that each ``torso'' has a (list-)coloring with bounded weak diameter, and finally show that the coloring of the ``torsos'' can be combined to obtain a list-coloring of the whole graph with bounded weak diameter.
Even though this strategy was also used in earlier work \cite{bbeglps,l_asdim} when proving the existence of a non-list-coloring with bounded weak diameter, there are extra challenges when using this strategy to prove Theorems \ref{minor_intro} and \ref{tw_intro} about list-coloring.

First, we have to develop a machinery for handling tree-decompositions to reduce the list-coloring problem for the whole graph to the one for the torsos as mentioned in the third step above.
This can be done by modifying the machinery for the non-list-coloring version by considering more complicated gadgets when applying induction.
However, as we will explain below, this machinery, which does not involve precoloring, is not strong enough in our application.

The main challenges in this paper arises from the aforementioned second step about finding list-colorings of the torsos.
This step is simpler in the previous work about non-list-coloring.
According to structure theorems, each torso ``essentially'' has bounded layered tree-width, and a non-list-3-coloring with bounded weak diameter for such a graph is not very hard to find: because the colors allowed for each vertex are the same, one can sacrifice a color for each layer to make sure that every monochromatic component of any coloring must be contained in the union of two consecutive layers, which induces a subgraph with bounded tree-width.
So the non-list-3-coloring problem is essentially reduced to the 2-coloring problem for graphs with bounded tree-width, which can be easily solved by the aforementioned machinery for handling tree-decompositions.

However, this trick does not work for list-coloring, because the lists for different vertices are very different.
In fact, we do not know how to prove that graphs with bounded layered tree-width are 3-choosable with bounded weak diameter, so the aforementioned structure theorem for $H$-minor free graphs about layered tree-width is too weak for us.
The structure theorem for $H$-minor free graphs we will use in this paper is the one proved by Robertson and Seymour \cite{rs_XVI}, where the torsos are ``near embeddings'' in surfaces of bounded Euler genus.
Here a near embedding consists of a surface part, vortices, and other additional elements that we do not discuss here as they are easy to be handled for weak diameter coloring.

To find a list-coloring of near embeddings with bounded weak diameter, we will use the result of Dvo\v{r}\'{a}k and Norin about list-coloring of graphs embedded in surfaces (Theorem \ref{surface_ch}) to color the surface part of a near embedding and then extend the coloring to the vortices.
To extend coloring from the surface part to the vortices, it suffices to consider the vortices together with the monochromatic components of the coloring of the surface part intersecting the vortices.
Note that the union of the vortices and those monochromatic components is a graph with bounded tree-width with a precoloring with bounded weak diameter. 

Therefore, the key challenge lies in proving Theorem \ref{tw_intro} about extending a precoloring with bounded weak diameter to a coloring with bounded weak diameter for $H$-minor free graphs with planar $H$.
In order to prove it, we develop a more technical machinery for handling tree-decompositions to reduce the list-coloring problem for the whole graph to the one for the torsos so that we can incorporate the precoloring with bounded weak diameter. 

In Section \ref{sec:tree}, we will prove this machinery for handling list-coloring problems with a given tree-decomposition, as long as the list-assignment of a graph $G$ is ``legitimate''.
Note that the list-assignments considered in our proof allow some vertices having list size 1.
Those vertices are exactly the precolored vertices.
``Legitimate'' list-assignments will be defined in Section \ref{sec:preparation}, which roughly says that the unique coloring for the vertices with list size 1 has bounded weak diameter.
Such a requirement for ``legitimate'' list-assignments is required because it is impossible to extend a precoloring with unbounded weak diameter to a coloring with bounded weak diameter.

In Section \ref{sec:main_color}, we show how to use the machinery developed in Section \ref{sec:tree} to prove Theorem \ref{tw_intro}, then show how to use Theorem \ref{tw_intro} to find a list-coloring for near embeddings, and finally prove Theorem \ref{minor_intro}.

In Section \ref{sec:odd_minors}, we consider problems about odd minors and prove Theorems \ref{bipartite_intro} and \ref{odd_minor_intro}. 
Finally we include some concluding remarks in Section \ref{sec:concluding}.

\section{Preparation} \label{sec:preparation}

\subsection{Notation}
Let $G$ be a graph.
Let $S$ be a subset of $V(G)$.
Then $G[S]$ is defined to be the subgraph of $G$ induced on $S$.
We define $N_G(S)$ to be the set of vertices not in $S$ but adjacent in $G$ to some vertices in $S$.
For a vertex $v$ of $G$, we define $N_G(v)$ to be $N_G(\{v\})$ and call every element in $N_G(v)$ a {\it neighbor} of $v$ in $G$.

\subsection{Centered sets}
Let $G$ be a graph.
For a set $S\subseteq V(G)$ and a nonnegative integer $r$, we define $N_G^{\leq r}[S] = \{v \in V(G): \dist_G(v,s) \leq r$ for some $s \in S\}$.
For any nonnegative integers $k$ and $r$, we say that a subset $Z$ of $V(G)$ is {\it $(k,r)$-centered in $G$} if $Z \subseteq N_G^{\leq r}[S]$ for some subset $S$ of $V(G)$ with $|S| \leq k$.

It is ``trivial'' to find a coloring with bounded weak diameter for a $(k,r)$-centered set, as observed in \cite{bbeglps}.

\begin{lemma}[{\cite[Observation 2.5]{bbeglps}}] \label{all_centered}
For any integers $k \geq 0,r \geq 0$, there exists an integer $N^*\ge 1$ such that for every graph $G$ and every integer $m\ge 1$, if $V(G)$ is $(k,r)$-centered, then any $m$-coloring of $G$ has weak diameter in $G$ at most $N^*$. 
\end{lemma} 

For $i \in \{1,2\}$, let $f_i$ be a function with domain $S_i$.
If $f_1(x)=f_2(x)$ for every $x \in S_1 \cap S_2$, then we define $f_1 \cup f_2$ to be the function with domain $S_1 \cup S_2$ such that $(f_1 \cup f_2)(x) = f_i(x)$ for any $i\in \{1,2\}$ and $x \in S_i$.

We will frequently use a special case of a result in \cite{l_ANdim}.

\begin{lemma}[{\cite[Special case of Lemma 3.6]{l_ANdim}}] \label{add_centered}
For any integers $k \geq 0,r \geq 0, N \geq 1$, there exists an integer $N^* \geq N$ such that the following holds.
Let $G$ be a graph.
Let $Z$ be a $(k,r)$-centered subset of $V(G)$ in $G$.
Let $m$ be a positive integer.
Let $R \subseteq V(G)$.
Let $c_Z: Z \rightarrow [m]$.
Let $c$ be an $m$-coloring of $G-(R \cup Z)$ with weak diameter in $G$ at most $N$.
Then the $m$-coloring $c \cup c_Z|_{Z-R}$ of $G-R$ has weak diameter in $G$ at most $N^*$.
\end{lemma}

\subsection{Legitimate list-assignments}

Let $(T,\X)$ be a tree-decomposition of a graph $G$, where $\X=(X_t: t \in X_t)$.
For every $S \subseteq V(T)$, we define $$X_S=\bigcup_{t \in S}X_t.$$
Similarly, for every subgraph $S$ of $T$, we define $X_S = X_{V(S)}$.
For every $S \subseteq V(G)$, we define $$\X \cap S = (X_t \cap S: t \in V(T));$$ note that $(T,\X \cap S)$ is a tree-decomposition of $G[S]$.

For any list-assignment $L$ of a graph $G$, we define $$1_L= \{v \in V(G): \lvert L(v) \rvert=1\}.$$

Let $G$ be a graph.
Let $(T,\X)$ be a tree-decomposition of $G$, where $\X$ is denoted by $(X_t: t \in V(T))$.
Let $m,s,r,k$ be integers.
We say that a list-assignment $L$ of $G$ is {\it $(T,\X,m,s,r,k)$-legitimate} if $L$ satisfies the following conditions:
    \begin{itemize}
        \item[(L1)] $|L(v)| \in \{1,m\}$ for every $v \in V(G)$.
        \item[(L2)] For every $t \in V(T)$, $1_L \cap X_t$ is $(s,r)$-centered in $G[X_t]$.
        \item[(L3)] The unique $L|_{1_L}$-coloring of $G[1_L]$ has weak diameter in $G$ at most $k$.
    \end{itemize}

\begin{lemma} \label{comp_legitimate}
Let $G$ be a graph.
Let $L$ be a $(T,\X,m,s,r,k)$-legitimate list-assignment of $G$ for some tree-decomposition $(T,\X)$ and integers $m,s,r,k$.
Let $C$ be a component of $G$.
Then $L|_{V(C)}$ is a $(T,\X \cap V(C),m,s,r,k)$-legitimate list-assignment of $C$.
\end{lemma}

\begin{pf}
Clearly $L|_{V(C)}$ satisfies (L1).

Let $t \in V(T)$.
Since $1_L \cap X_t$ is an $(s,r)$-centered set in $G[X_t]$, there exists $S_t \subseteq X_t$ with $|S_t| \leq s$ such that $1_L \cap X_t \subseteq N_{G[X_t]}^{\leq r}[S_t]$.
Then $1_{L|_{V(C)}} \cap X_t \cap V(C) = 1_L \cap X_t \cap V(C) \subseteq N_{G[X_t]}^{\leq r}[S_t] \cap V(C) \subseteq N_{G[X_t \cap V(C)]}^{\leq r}[S_t \cap V(C)] = N_{C[X_t \cap V(C)]}^{\leq r}[S_t \cap V(C)]$ is $(s,r)$-centered in $C[X_t \cap V(C)]$.
Hence (L2) holds.

Let $L'=L|_{V(C)}$.
Since every monochromatic component of the unique $L'|_{1_{L'}}$-coloring of $C[1_{L'}]$ is a monochromatic component of the unique $L|_{1_L}$-coloring of $G[1_L]$ contained in $C$, (L3) holds.
\end{pf}

\subsection{Extendable list-assignments}

A {\it rooted tree} is a directed graph whose underlying graph is a tree, where exactly one vertex, called the \emph{root}, has in-degree 0, and all other vertices have in-degree $1$.
A {\it rooted tree-decomposition} of a graph $G$ is a tree-decomposition $(T,\X)$ of $G$ such that $T$ is a rooted tree.

Let $(T,\X)$ be a rooted tree-decomposition of a graph $G$ with $\X=(X_t: t \in V(T))$.
For $t \in V(T)$, a {\it child-extension at $t$ with respect to $(T,\X)$} is a graph obtained from $G[X_t]$ by, for each child $t'$ of $t$ in $T$, adding new vertices and adding new edges that are either between new vertices or between a new vertex and a vertex in $X_t \cap X_{t'}$ such that every component of the subgraph induced by the new vertices contains at most two vertices.

Let $G$ be a graph.
Let $(T,\X)$ be a rooted tree-decomposition of $G$ with $\X=(X_t: t \in V(T))$.
Let $m,s,r,k,N$ be integers.
We say that a list-assignment $L$ of $G$ is {\it $(T,\X,m,s,r,k, \allowbreak N)$-extendable} if the following conditions hold:
    \begin{itemize}
        \item[(E1)] $L$ is $(T,\X,m,s,r,k)$-legitimate.
        \item[(E2)] For any $t \in V(T)$, any child-extension $G_t$ at $t$ with respect to $(T,\X)$, and any list-assignment $L_t$ of $G_t$ with $L_t|_{X_t}=L|_{X_t}$, there exists an $L_t$-coloring of $G_t$ with weak diameter in $G_t$ at most $N$.
    \end{itemize}

Let $G$ be a graph.
Let $(T,\X)$ be a rooted tree-decomposition of $G$.
Let $m,s,r,k,N$ be integers.
We say that a list-assignment $L$ of $G$ is {\it hereditarily $(T,\X,m,s,r,k,N)$-extendable} if the following conditions hold:
    \begin{itemize}
        \item[(H1)] $L$ is $(T,\X,m,s,r,k)$-legitimate.
        \item[(H2)] For every induced subgraph $H$ of $G$ and for every list-assignment $L_H$ of $H$ with $L_H(v) \supseteq L(v)$ for every $v \in V(H)$, if $L_H$ is $(T,\X \cap V(H),m,s,r,k)$-legitimate, then $L_H$ is $(T,\X \cap V(H),m,s,r,k,N)$-extendable. 
    \end{itemize}

\begin{lemma} \label{hereditary_list}
Let $(T,\X)$ be a rooted tree-decomposition of a graph $G$, where $\X=(X_t: t \in V(T))$.
Let $T'$ be a subtree of $T$, and let $\X'=(X_t: t \in V(T'))$. 
Let $G' = G[\bigcup_{t \in V(T')}X_{t}]$.
Let $L$ be a hereditarily $(T,\X,m,s,r,k,N)$-extendable list-assignment of $G$ for some integers $m,s,r,k,N$. 
Let $L'$ be a list-assignment of $G'$ such that $L'(v) \supseteq L(v)$ for every $v \in V(G')$, and $1_{L'} \cap X_t=\emptyset$ for every $t \in V(T)-V(T')$.
If $L'$ is $(T',\X',m,s,r,k)$-legitimate, then $L'$ is hereditarily $(T',\X',m,s,r,k,N)$-extendable.
\end{lemma}

\begin{pf}
Since $L'$ is $(T',\X',m,s,r,k)$-legitimate, it suffices to verify (H2) for $L'$.
Let $H$ be an induced subgraph of $G'$.
Let $L_H$ be a $(T',\X' \cap V(H),m,s,r,k)$-legitimate list-assignment of $H$ such that $L_H(v) \supseteq L'(v)$ for every $v \in V(H)$.
It suffices to show that $L_H$ is $(T',\X' \cap V(H),m,s,r,k,N)$-extendable.

Let $t' \in V(T')$.
Let $H_{t'}$ be a child-extension at $t'$ with respect to $(T', \X' \cap V(H))$.
Let $L_{t'}$ be a list-assignment of $H_{t'}$ with $L_{t'}|_{X_{t'} \cap V(H)} = L_H|_{X_{t'} \cap V(H)}$.
It suffices to show that there exists an $L_{t'}$-coloring of $H_{t'}$ with weak diameter in $H_{t'}$ at most $N$.

We first show that $L_H$ is $(T,\X \cap V(H),m,s,r,k)$-legitimate.
Since $L_H$ is $(T',\X' \cap V(H),m,s,r,k)$-legitimate, (L1) and (L3) hold for $L_H$ being $(T,\X \cap V(H),m,s,r,k)$-legitimate, and for every $t \in V(T')$, $1_{L_H} \cap X_t \cap V(H)$ is $(s,r)$-centered in $H[X_t \cap V(H)]$.
And for every $t \in V(T)-V(T')$, $1_{L_H} \cap X_t \cap V(H) \subseteq 1_{L'} \cap X_t = \emptyset$ is $(s,r)$-centered in $H[X_t \cap V(H)]$.
Hence $L_H$ is $(T,\X \cap V(H),m,s,r,k)$-legitimate.

Note that $H$ is also an induced subgraph of $G$, and $L_H(v) \supseteq L'(v) \supseteq L(v)$ for every $v \in V(H)$.
Since $L_H$ is $(T,\X \cap V(H),m,s,r,k)$-legitimate, by (H2) for $L$ being hereditarily $(T,\X,m,s,r,k,N)$-extendable, $L_H$ is $(T,\X \cap V(H),m,s,r,k,N)$-extendable.

Since $t' \in V(T') \subseteq V(T)$ and $T' \subseteq T$, we know that $H_{t'}$ is also and child-extension at $t'$ with respect to $(T,\X \cap V(H))$.
Since $L_H$ is $(T,\X \cap V(H),m,s,r,k,N)$-extendable and $L_{t'}|_{X_{t'} \cap V(H)} = L_H|_{X_{t'} \cap V(H)}$, by (E2), we know that there exists an $L_{t'}$-coloring of $H_{t'}$ with weak diameter in $H_{t'}$ at most $N$.
This proves the lemma.
\end{pf}

\bigskip

Let $(T,\X)$ be a rooted tree-decomposition of a graph $G$, where $\X=(X_t: t \in V(T))$.
Let $e$ be an edge of $T$.
Let $T_e$ be the component of $T-e$ disjoint from the root of $T$.
The {\it truncation of $(T,\X)$ at $e$} is a pair $(T_e',\X')$ such that
	\begin{itemize}
		\item $T_e'$ is the rooted tree obtained from $T_e$ by adding a new node $t^*_e$ adjacent to the root of $T_e$, where $t^*_e$ is the root of $T_e'$, and
		\item $\X'=(X'_t: t \in V(T_e'))$, where 
			\begin{itemize}
				\item $X'_{t^*_e}=X_e$, where $X_e$ is the intersection of the bags at the ends of $e$ in $(T,\X)$, and
				\item $X'_t=X_t$ for every $t \in V(T_e)$.
			\end{itemize}
	\end{itemize}
Note that $G[\bigcup_{t \in V(T_e')}X'_{t}] = G[\bigcup_{t \in V(T_e)}X_{t}]$, and $(T_e',\X')$ is its tree-decomposition.

\begin{lemma} \label{hereditary_list_2}
Let $(T,\X)$ be a rooted tree-decomposition of a graph $G$, where $\X=(X_t: t \in V(T))$.
Let $e$ be an edge of $T$, and let $(T_e',\X')$ be the truncation of $(T,\X)$ at $e$.
Let $G_e = G[\bigcup_{t \in V(T_e')}X'_{t}]$.
Let $L$ be a hereditarily $(T,\X,m,s,r,k,N)$-extendable list-assignment of $G$ for some integers $m,s,r,k,N$. 
Let $L'$ be a list-assignment of $G_e$ such that $L'(v) \supseteq L(v)$ for every $v \in V(G_e)$, and $1_{L'} \cap X_t=\emptyset$ for every $t \in V(T)-V(T_e')$.
If $L'$ is $(T_e',\X',m,s,r,k)$-legitimate, then $L'$ is hereditarily $(T_e',\X',m,s,r,k,N)$-extendable.
\end{lemma}

\begin{pf}
Since $L'$ is $(T_e',\X',m,s,r,k)$-legitimate, to show that $L'$ is hereditarily $(T_e',\X',m,s,r, \allowbreak k,N)$-extendable, it suffices to verify (H2).
Let $H$ be an induced subgraph of $G_e$.
Let $L_H$ be a $(T_e',\X' \cap V(H),m,s,r,k)$-legitimate list-assignment of $H$ such that $L_H(v) \supseteq L'(v)$ for every $v \in V(H)$.
It suffices to show that $L_H$ satisfies (E2) for being $(T_e',\X' \cap V(H),m,s,r,k,N)$-extendable.
Let $t \in V(T_e')$.
Let $G_t$ be a child-extension at $t$ with respect to $(T_e',\X' \cap V(H))$.
Let $L_t$ be a list-assignment of $G_t$ with $L_t|_{X'_t \cap V(H)}=L_H|_{X'_t \cap V(H)}$.
It suffices to show that there exists an $L_t$-coloring of $G_t$ with weak diameter in $G_t$ at most $N$.

Let $t^*_e$ be the root of $T_e'$.
Let $T_e = T_e'-t^*_e$.
Note that $(T_e,(X'_t: t \in V(T_e)))$ is a tree-decomposition of $G_e$.

Since $L'$ is $(T_e',\X',m,s,r,k)$-legitimate and $(T_e, (X'_t: t \in V(T_e)))$ can be obtained from $(T_e',\X')$ by removing the root and its bag, we know that $L'$ is $(T_e,(X'_t: t \in V(T_e)),m,s,r,k)$-legitimate.
By Lemma \ref{hereditary_list} (taking $T'=T_e$ and $G'=G_e$), 
	\begin{itemize}
		\item[(i)] $L'$ is hereditarily $(T_e,(X'_t: t \in V(T_e)),m,s,r,k,N)$-extendable.
	\end{itemize}

If $t \in V(T_e)$, then $G_t$ is also a child-extension at $t$ with respect to $(T_e,(X'_z \cap V(H): z \in V(T_e)))$, so there exists an $L_t$-coloring of $G_t$ with weak diameter in $G_t$ at most $N$ by (i). 
So we may assume $t = t^*_e$.

Let $p$ be the end of $e$ not in $T_e$.
Since $t = t^*_e$ and $H \subseteq G_e$, we have $X'_t \cap V(H) = X_p \cap V(G_e) \cap V(H) = X_p \cap V(H)$, so 
	\begin{itemize}
		\item[(ii)] $G_t$ is a child-extension at $p$ with respect to $(T,\X \cap V(H))$.
	\end{itemize}

Since $H$ is an induced subgraph of $G_e$, $H$ is also an induced subgraph of $G$, and $L_H(v) \supseteq L'(v) \supseteq L(v)$ for every $v \in V(H)$.
Since $L_H$ is $(T_e',\X' \cap V(H),m,s,r,k)$-legitimate, $L_H$ satisfies (L1) and (L3) for being $(T,\X \cap V(H), m,s,r,k)$-legitimate.
Since $1_{L_H} \subseteq 1_{L'}$ and $1_{L'} \cap X_t=\emptyset$ for every $t \in V(T)-V(T_e)$, $L_H$ satisfies (L2) for being $(T,\X \cap V(H), m,s,r,k)$-legitimate.
Hence $L_H$ is $(T,\X \cap V(H), m,s,r,k)$-legitimate.
Since $L$ is hereditarily $(T,\X,m,s,r,k,N)$-extendable, $L_H$ is $(T,\X \cap V(H),m,s,r,k,N)$-extendable.

Recall that $X_p \cap V(H) = X'_t \cap V(H)$.
Since $L_t|_{X'_t \cap V(H)}=L_H|_{X'_t \cap V(H)}$, we have $L_t|_{X_p \cap V(H)} =L_H|_{X_p \cap V(H)}$.
Recall that $G_t$ is a child-extension at $p$ with respect to $(T,\X \cap V(H))$ by (ii).
Since $L_t|_{X_p \cap V(H)}=L_H|_{X_p \cap V(H)}$, the condition (E2) for $L_H$ being $(T,\X \cap V(H),m,s,r,k,N)$-extendable implies that there exists an $L_t$-coloring of $G_t$ with weak diameter in $G_t$ at most $N$.
This proves the lemma.
\end{pf}

\section{Coloring a graph with a construction} \label{sec:tree}

\subsection{Overview}
The goal of this section is to prove the main technical lemma (Lemma \ref{tree_extension_list}) of this paper.
The main purpose of this lemma is to show that if a graph $G$ has a tree-decomposition $(T,\X)$ of bounded adhesion, then $G$ has an $L$-coloring with bounded weak diameter for every hereditarily $(T,\X,m,s,r,k,N)$-extendable list-assignment $L$.
In order to have sufficiently strong inductive hypothesis, we consider a stronger statement that allows vertices that are not far away from the root bag to be precolored.
We sketch a strategy for proving this statement, and the utility of precoloring vertices will be clear.

Let $Z$ be the set consisting of precolored vertices.
For simplicity, we assume that $G$ is connected and $Z$ consists of all vertices that are not far away from the root bag in the proof sketch.
Consider the subgraph $T_0$ of $T$ induced by the tree nodes whose bags intersect $Z$.
Since $Z$ consists of all vertices that are not far away from the root bag, $T_0$ is a subtree of $T$ containing the root of $T$.
Let $G_0$ be the subgraph of $G$ induced by $X_{T_0}$.
Note that $G_0-Z$ has a tree-decomposition of adhesion smaller than $(T,\X)$.
So we can apply induction on the adhesion of $(T,\X)$ to obtain a coloring of $G_0-Z$ with bounded weak diameter, and then combining this coloring with the precoloring on $Z$ to obtain a coloring $c_0$ of $G_0$ with bounded weak diameter in $G_0$ by using Lemma \ref{add_centered}.

Let $T_1 = T-V(T_0)$, and let $G_1$ be the subgraph of $G$ induced by $X_{T_1}$.
For every component $C$ of $T_1$, the coloring $c_0$ gives a precoloring on $X_C \cap V(G_0)$, and $X_C \cap V(G_0)$ is contained in the root bag of $(C,\X \cap X_C)$; since $|X_C \cap V(G_0)|$ is bounded, we can further precolor vertices in $X_C$ that are not far away from $X_C \cap V(G_0)$ such that for every coloring of $G_1$ obtained by extending this precoloring, every monochromatic component intersecting $X_C \cap V(G_0)$ is contained in the set of precolored vertices; note that the currently precolored vertices are not far away from the root bag of $(C, \X \cap X_C)$, so we can apply induction on $|V(G)|$ to obtain a coloring $c_C$ of $G[X_C]$ with bounded weak diameter.
By combining the coloring $c_C$ for each component $C$ of $T_1$, we obtain a coloring $c_1$ of $G_1$ with bounded weak diameter such that the union of the monochromatic components intersecting $V(G_0) \cap V(G_1)$ is not far away from $V(G_0) \cap V(G_1)$.

Then we hope that the coloring $c_0 \cup c_1$ has bounded weak diameter in $G$ to finish the proof.
However, the weak diameter in $G$ of $c_0 \cup c_1$ can be arbitrarily large because some single $(c_0 \cup c_1)$-monochromatic component possibly contains an unbounded number of $c_0$-monochromatic components. 
To resolve this issue, we construct a graph $G_0'$ by attaching ``gadgets'' to $G_0$ on $V(G_0) \cap V(G_1)$ to simulate paths in $G_1$ between vertices in $V(G_0) \cap V(G_1)$ without ``significantly changing'' the original distance in $G$ between any two vertices in $V(G_0) \cap V(G_1)$; the added gadgets ensure that if we obtain a coloring $c_0'$ of $G_0'$ with bounded weak diameter in $G_0'$, then we can precolor the vertices in $G_1$ not far away from $V(G_0) \cap V(G_1)$ according to the colors of the gadgets in $c_0'$ such that the coloring obtained by combining the restriction of $c_0'$ on $G_0$ with the resulting coloring $c_1$ of $G_1$ obtained by extending the precoloring has bounded weak diameter in $G$.

Note that we can construct a tree-decomposition of $G_0'$ by adding leave bags to $(T,\X)$ to accommodate those gadgets without increasing the adhesion.
But the size of the neighborhood of a gadget can be as large as the adhesion of $(T,\X)$, and $Z$ is disjoint from such a neighborhood.
So deleting $Z$ from $G_0'$ does not decrease the adhesion of the corresponding tree-decomposition of $G_0'$.
Hence we cannot apply induction on adhesion in a naive way as we outline above.
To remedy this, we consider $(\eta,\theta)$-constructions defined in Section \ref{subsec:proof_main_lemma}, which are tree-decompositions of adhesion at most $\theta \geq \eta$ but only special adhesion sets have size greater than $\eta$.

Lemma \ref{tree_extension_list} is the formal form of the key lemma stated in terms of $(\eta,\theta)$-constructions.
The proof of Lemma \ref{tree_extension_list} follows from the strategy sketched above.

\subsection{Proof} \label{subsec:proof_main_lemma}

Let $\eta,\theta$ be nonnegative integers with $\eta \leq \theta$.
Let $G$ be a graph.
An {\it $(\eta,\theta)$-construction} of $G$ is a rooted tree-decomposition $(T,\X)$ of $G$ with adhesion at most $\theta$ satisfying the following conditions:
	\begin{itemize}
		\item[(C1)] For every edge $tt' \in E(T)$, if $\lvert X_t \cap X_{t'} \rvert > \eta$, then one end of $tt'$ has no child, say $t'$, and every component of $G[X_{t'}-X_t]$ has at most two vertices. 
		\item[(C2)] For the root $t^*$ of $T$, 
		\begin{itemize}
		    \item $\lvert X_{t^*} \rvert \leq \theta$, and
			\item if $\eta>0$, then $X_{t^*} \neq \emptyset$, 
		\end{itemize}
	\end{itemize}
We say that $G$ is {\it $(\eta,\theta)$-constructible} if there exists an $(\eta,\theta)$-construction of $G$.

\begin{lemma} \label{tree_extension_list}
For any integers $\theta \geq 0,s \geq \theta, r \geq 1, k \geq 1$, there exists a function $f^*: ({\mathbb N} \cup \{0\}) \times {\mathbb N} \rightarrow {\mathbb N}$ such that the following holds.
Let $m \geq 2$ and $N \geq 4$ be integers.
Let $\eta$ be an integer with $0 \leq \eta \leq \theta$.
Let $G$ be an $(\eta,\theta)$-constructible graph with an $(\eta,\theta)$-construction $(T,\X)$. 
Denote $\X$ by $(X_t: t \in V(T))$.
Let $t^*$ be the root of $T$.
If $L$ is a hereditarily $(T,\X,m,s,r,k,N)$-extendable list-assignment of $G$, then for every $Z \subseteq N_G^{\leq (k+2)(\theta+1)}[X_{t^*}]$, every $L|_{Z}$-coloring $c_Z$ of $G[Z]$ can be extended to an $L$-coloring of $G$ with weak diameter in $G$ at most $f^*(\eta,N)$.
\end{lemma}

\begin{pf}
Let $\theta \geq 0, s \geq \theta, r \geq 1, k \geq 1$ be integers.
We define the following.
	\begin{itemize}
        \item Let $N_1$ be the number $N^*$ mentioned in Lemma \ref{all_centered} by taking $(k,r)=(\theta,(k+2)(\theta+1))$.
        \item Let $f_1: {\mathbb N} \rightarrow {\mathbb N}$ be the function such that for every $x \in {\mathbb N}$, $f_1(x)$ is the integer $N^*$ mentioned in Lemma \ref{add_centered} by taking $(k,r,N)=(\theta,(k+2)(\theta + 1),x)$. 
        \item Let $f_2: {\mathbb N} \rightarrow {\mathbb N}$ be the function such that for every $x \in {\mathbb N}$, $f_2(x)$ is the integer $N^*$ mentioned in Lemma \ref{add_centered} by taking $(k,r,N)=(s,r,x)$.
        \item Let $f_3: {\mathbb N} \rightarrow {\mathbb N}$ be the function such that for every $x \in {\mathbb N}$, $f_3(x)$ is the integer $N^*$ mentioned in Lemma \ref{add_centered} by taking $(k,r,N)=(\theta,(k+2)(\theta + 1)+k+r,x)$. 
		\item Define $f^*: ({\mathbb N} \cup \{0\}) \times {\mathbb N} \rightarrow {\mathbb N}$ to be the function such that for every $y \in {\mathbb N}$,
		\begin{itemize}
			\item $f^*(0,y)= N_1+f_1(y)$, and 
			\item for any $x \in {\mathbb N}$, $f^*(x,y)= (k+2)(\theta+1)^2 \cdot (4+f_3(f_1(f^*(x-1,N_1+f_2(y)))))$.  
		\end{itemize}
	\end{itemize}
Note that $f_i(x) \geq x$ for any $i \in [3]$ and $x \in {\mathbb N}$, so $f^*(x,y)$ is increasing on $x$.
	
Let $m,N,\eta,G,(T,\X),t^*,L,Z,c_Z$ be as defined in the lemma.
Suppose to the contrary that $c_Z$ cannot be extended to an $L$-coloring of $G$ with weak diameter in $G$ at most $f^*(\eta,N)$.
That is, $m,N,\eta,G,(T,\X),t^*,L,Z,c_Z$ provide a counterexample to this lemma.
We further assume that $m,N,\eta,G,(T,\X),t^*,L,Z,c_Z$ are chosen so that the lexicographical order of the tuple $(\eta,|V(G)-Z|+|V(G)|,|V(T)|)$ is as small as possible. 

Since $(T,\X)$ is an $(\eta,\theta)$-construction, $|X_{t^*}| \leq \theta$.
So $Z$ is $(\theta,(k+2)(\theta+1))$-centered.
If $Z=V(G)$, then $c_Z$ is itself an $L$-coloring of $G$ with weak diameter in $G$ at most $N_1 \leq f^*(\eta,N)$ by Lemma \ref{all_centered}, a contradiction. 
So $Z \neq V(G)$. 
	
\begin{clm} \label{cl:1} 
$\eta \geq 1$. 
\end{clm}
	
\begin{claimproof}
Suppose to the contrary that $\eta=0$. 
If for every component $C$ of $G$, $c_{Z}|_{Z \cap V(C)}$ can be extended to an $L|_{V(C)}$-coloring of $C$ with weak diameter in $G$ at most $f^*(\eta,N)$, then $c_Z$ can be extended to an $L$-coloring of $G$ with weak diameter in $G$ at most $f^*(\eta,N)$, a contradiction.
So there exists a component $C$ of $G$ such that $c_{Z}|_{Z \cap V(C)}$ cannot be extended to an $L|_{V(C)}$-coloring of $C$ with weak diameter in $G$ at most $f^*(\eta,N)$.
  
Let $W=\{tt' \in E(T): X_t \cap X_{t'}=\emptyset\}$. 
Then $V(C) \subseteq X_{T_C}$ for some component $T_C$ of $T-W$. 
Let $t_C$ be the root of $T_C$.
Let $G_C=G[X_{T_C}]$.
For every edge $tt'\in E(T_C)$, $|X_t\cap X_{t'}|>0 = \eta$, so by (C1), one end of $tt'$ has no child, say $t'$, and every component of $G[X_{t'}-X_t]$ contains at most two vertices.  
In particular, $T_C$ is a star, and $G_C$ is a child-extension at $t_C$ with respect to $(T,\X)$. 

Since $L$ is hereditarily $(T,\X,m,s,r,k,N)$-extendable, by taking $H=G$ and $L_H=L$ in (H2), (H1) implies that $L$ is $(T,\X,m,s,r,k,N)$-extendable.
By (E2), since $L|_{V(G_C)}$ is a list-assignment of $G_C$ with $(L|_{V(G_C)})|_{X_{t_C}}=L|_{X_{t_C}}$, there exists an $L|_{V(G_C)}$-coloring $c_C$ of $G_C$ with weak diameter in $G_C$ at most $N$.
Since $G_C \subseteq G$, $c_C$ has weak diameter in $G$ at most $N$.
Since every $c_C|_{V(C)-Z}$-monochromatic component is contained in a $c_C$-monochromatic component, $c_C|_{V(C)-Z}$ has weak diameter in $G$ at most $N$.
Since $Z$ is $(\theta,(k+2)(\theta+1))$-centered in $G$, so is $V(C) \cap Z$.
So by Lemma \ref{add_centered} (taking $(k,r,N,G,Z,m,R,c_Z,c)=(\theta,(k+2)(\theta+1),N,G,V(C) \cap Z, |\bigcup_{v \in V(G)}L(v)|, V(G)-V(C), c_Z|_{V(C) \cap Z},c_C|_{V(C)-Z})$), we know that $c_C|_{V(C)-Z} \cup c_Z|_{V(C) \cap Z}$ is an $L|_{V(C)}$-coloring of $C$ with weak diameter in $G$ at most $f_1(N) \leq f^*(\eta,N)$, a contradiction.
\end{claimproof}

\begin{clm} \label{clm:conn}
$G$ is connected. 
\end{clm}
	
\begin{claimproof}
Suppose to the contrary that $G$ is disconnected.
Note that there exists a component $C$ of $G$ such that $c_Z|_{V(C) \cap Z}$ cannot be extended to an $L|_{V(C)}$-coloring of $C$ with weak diameter in $G$ (and hence in $C$) at most $f^*(\eta,N)$. 
Since $G$ is disconnected, $|V(C)|<|V(G)|$.
		
Since $L$ is hereditarily $(T,\X,m,s,r,k,N)$-extendable, (H1) implies that $L$ is $(T,\X,m,s,r, \allowbreak k)$-legitimate.
By Lemma \ref{comp_legitimate}, $L|_{V(C)}$ is $(T,\X \cap V(C),m,s,r,k)$-legitimate.
Since every induced subgraph of $C$ is an induced subgraph of $G$ and $L$ is hereditarily $(T,\X,m,s,r,k,N)$-extendable, (H2) implies that $L|_{V(C)}$ is hereditarily $(T,\X \cap V(C),m,s,r,k,N)$-extendable.

Let $T_C$ be the subtree of $T$ induced by $\{t \in V(T): X_t \cap V(C) \neq \emptyset\}$. 
Let $\X_C = (X_t \cap V(C): t \in V(T_C))$.
Since $(T_C,\X_C)$ can be obtained from $(T,\X \cap V(C))$ by removing nodes with empty bags, $L|_{V(C)}$ is $(T_C,\X_C,m,s,r,k)$-legitimate.
By Lemma \ref{hereditary_list} (taking $(G,T,\X,T',\X',G',L,L')=(C,T,\X \cap V(C),T_C,\X_C,C,L|_{V(C)},L|_{V(C)})$), $L|_{V(C)}$ is hereditarily $(T_C,\X_C,m,s,r,k,N)$-extendable.

Suppose $t^* \in V(T_C)$.
Then $t^*$ is the root of $T_C$ and $X_{t^*} \cap V(C) \neq \emptyset$.
So $(T_C,\X_C)$ is an $(\eta,\theta)$-construction of $C$ with $|V(C)-(Z \cap V(C))|+|V(C)|< |V(G)-Z|+|V(G)|$.
Note that $Z \cap V(C) \subseteq N_{G}^{\leq (k+2)(\theta+1)}[X_{t^*}] \cap V(C) \subseteq N_C^{\leq (k+2)(\theta+1)}[X_{t^*} \cap V(C)]$.
By the minimality of $(\eta,|V(G)-Z|+|V(G)|)$, we know that $c_Z|_{V(C) \cap Z}$ can be extended to an $L|_{V(C)}$-coloring of $C$ with weak diameter in $C$ at most $f^*(\eta,N)$, a contradiction.

So $t^* \not \in V(T_C)$. 
Let $T_C'$ be the rooted tree obtained from $T_C$ by adding a node $t^*_C$ adjacent to the root of $T_C$, where $t^*_C$ is the root of $T_C'$.
For every $t \in V(T_C)$, let $X'_t=X_t$.
Let $X'_{t^*_C}$ be a set consisting of a single vertex in the intersection of $V(C)$ and the bag of the root of $T_C$.
Let $\X'_C = (X_t': t \in V(T_C'))$.
By Claim \ref{cl:1}, $\theta \geq \eta \geq 1 = |X'_{t^*_C}|$, so $(T_C', \X_C')$ is an $(\eta,\theta)$-construction of $C$.

We shall show that $L|_{V(C)}$ is hereditarily $(T_C',\X_C',m,s,r,k,N)$-extendable.

Since $L|_{V(C)}$ is hereditarily $(T_C,\X_C,m,s,r,k,N)$-extendable, $L|_{V(C)}$ is $(T_C,\X_C,m,s,r,k)$-legitimate.
By Claim 1, $|X'_{t^*_C}| \leq 1 \leq \theta \leq s$, so $1_{L|_{V(C)}} \cap X'_{t^*_C}$ is $(s,r)$-centered in $C[X'_{t^*_C}]$.
Hence $L|_{V(C)}$ is $(T_C',\X_C',m,s,r,k)$-legitimate.

Note that for every induced subgraph $H$ of $C$, every child-extension at a node $t \in V(T_C)$ with respect to $(T_C',\X_C' \cap V(H))$ is also a child-extension at $t$ with respect to $(T_C,\X_C \cap V(H))$; for every induced subgraph $H$ of $C$, since $|X'_{t^*_C} \cap V(H)| \leq 1$, every child-extension at $t_C^*$ with respect to $(T_C',\X_C' \cap V(H))$ is a graph whose every component has radius at most 2, so every its coloring has weak diameter at most $4 \leq N$. 
Hence $L|_{V(C)}$ is hereditarily $(T_C',\X_C',m,s,r,k,N)$-extendable.

Recall $\lvert V(C) \rvert < \lvert V(G) \rvert$.
So $|V(C)-(Z \cap V(C))|+|V(C)| < |V(G)-Z|+|V(G)|$. 
By the minimality of $(\eta,|V(G)-Z|+|V(G)|)$, $c_Z|_{V(C) \cap Z}$ can be extended to an $L|_{V(C)}$-coloring of $C$ with weak diameter in $C$ at most $f^*(\eta,N)$, a contradiction.
\end{claimproof}
	
\begin{clm} \label{claim_Zbasic}
$X_{t^*} \neq \emptyset$, $Z=N_G^{\leq (k+2)(\theta+1)}[X_{t^*}]$ and $Z-X_{t^*} \neq \emptyset$. 
\end{clm}
	
\begin{claimproof}
By Claim 1 and (C2), $X_{t^*} \neq \emptyset$.
If there exists $v \in N_G^{\leq (k+2)(\theta+1)}[X_{t^*}]-Z$, then let $Z' = Z \cup \{v\}$ and let $c'$ be the function with domain $Z'$ obtained from $c_Z$ by further defining $c'(v)$ to be an arbitrary element in $L(v)$. 
Note that $|V(G)-Z'|<|V(G)-Z|$, so by the minimality of $(\eta,|V(G)-Z|+|V(G)|)$, we know that $c'$ (and hence $c_Z$) can be extended to an $L$-coloring of $G$ with weak diameter in $G$ at most $f^*(\eta,N)$, a contradiction.
Hence $Z=N_G^{\leq (k+2)(\theta+1)}[X_{t^*}]$. 
In particular, since $Z\neq V(G)$ and $G$ is connected by Claim \ref{clm:conn}, we find $Z-X_{t^*}\neq \emptyset$. 
\end{claimproof}

For each $e=tt' \in E(T)$, define $X_e=X_t \cap X_{t'}$. 
Since the tree-decomposition $(T,\X)$ has adhesion at most $\theta$, we know $|X_e|\le \theta$ for every $e\in E(T)$.
For each $e \in E(T)$, let $T_e$ be the component of $T-e$ disjoint from $t^*$. 

\begin{clm} \label{clm:X_e_nonempty}
For every $e \in E(T)$, $X_e \neq \emptyset$.    
\end{clm}

\begin{claimproof}
Suppose to the contrary that there exists $e \in E(T)$ with $X_e =\emptyset$.
Since $G$ is connected (by Claim \ref{clm:conn}) and $X_{t^*} \neq \emptyset$ (by Claim 3), $X_{T_e}=\emptyset$.
Let $T'=T-V(T_e)$.
Let $\X'=(X_t: t \in V(T'))$.
Since $(T,\X)$ is an $(\eta,\theta)$-construction of $G$, $(T',\X')$ is an $(\eta,\theta)$-construction of $G$.
Since $L$ is hereditarily $(T,\X,m,s,r,k,N)$-extendable, $L$ is $(T,\X,m,s,r,k)$-legitimate, so $L$ is $(T',\X',m,s,r,k)$-legitimate.
By Lemma \ref{hereditary_list} (with taking $L'=L$), $L$ is hereditarily $(T',\X',m,s,r,k,N)$-extendable.
Then by the minimality of $(\eta,|V(G)-Z|+|V(G)|,|V(T)|)$, $c_Z$ can be extended to an $L$-coloring of $G$ with weak diameter in $G$ at most $f^*(\eta,N)$, a contradiction.
\end{claimproof}

For each $v \in X_{t^*}$, let $T_v$ be the subgraph of $T$ induced by $$\{t \in V(T): N_G^{\leq (k+2)(\theta+1)}[\{v\}] \cap X_t \neq \emptyset\}.$$
Since $(T,\X)$ is a tree-decomposition, $T_v$ is a subtree of $T$ containing $t^*$.
So $\bigcup_{v \in X_{t^*}}T_v$ is a subtree of $T$ containing $t^*$.
	
Let $T_0 = \bigcup_{v \in X_{t^*}}T_v$.
Since $Z = N_G^{\leq (k+2)(\theta+1)}[X_{t^*}]$, we know $Z \subseteq X_{T_0}$.
Let $$U_E = \{e \in E(T): \text{exactly one end of $e$ is in } V(T_0)\}.$$ 
Note that $Z \cap X_{T_e}=\emptyset$ for every $e \in U_E$. 
	
Let $T_1=\bigcup_{e\in U_E}T_e$.
Note that $V(T_1)=V(T)-V(T_0)$.
Define $G_0=G[X_{T_0}]$, $G_1=G[X_{T_1}]$, and $G_e=G[X_{T_e}]$ for every edge $e\in U_E$.  
By Claim~\ref{claim_Zbasic}, $Z \cap X_t \neq \emptyset$ for every $t \in V(T_0)$, and $Z \cap X_t = \emptyset$ for every $t \in V(T_1)$.	
	
For each $e \in U_E$, we define a partition $\P_e$ of $X_e$ such that two vertices $x,y\in X_e$ are in the same part of $\P_e$ if and only if there exists a sequence $a_1,\dots,a_\theta$ of (not necessarily distinct) elements of $X_e$ such that $a_1=x$, $a_\theta=y$, and for every $i \in [\theta-1]$, there exists a path in $G_e$ from $a_i$ to $a_{i+1}$ of length at most $2(k+2)(\theta+1)+1$.
	
Define $H$ to be the graph obtained from $G_0$ by
    \begin{itemize}
        \item for each $e \in U_E$ and each $Y \in \P_e$, adding ${|\bigcup_{v \in V(G)}L(v)| \choose m}$ new vertices $v_{Y,1},v_{Y,2},..., \allowbreak v_{Y,{|\bigcup_{v \in V(G)}L(v)| \choose m}}$ and new edges such that $N_H(v_{Y,i})=Y$ for every $i \in [{|\bigcup_{v \in V(G)}L(v)| \choose m}]$, and 
        \item for each $e \in U_E$ and each monochromatic component $M$ of the unique $L|_{1_L}$-coloring of $G[L|_{1_L}]$ intersecting $N_{G_e}^{\leq 1}[X_e]-X_e$, adding two new vertices $u_{e,M}$ and $u'_{e,M}$ and adding new edges such that 
            \begin{itemize}
                \item $N_H(u_{e,M}) = \{u'_{e,M}\} \cup (X_e \cap N_{G_e}^{\leq 1}[V(M) \cap X_{T_e}-X_e])$, and  
                \item $N_H(u'_{e,M}) = \{u_{e,M}\} \cup (X_e \cap N_{G_e}^{\leq k}[V(M) \cap X_{T_e}-X_e])$.
            \end{itemize}
    \end{itemize}
	
For every $e \in U_E$, we define 
    \begin{itemize}
        \item $S_{e,Y}=\{v_{Y,i}: i \in [{\lvert \bigcup_{v \in V(G)}L(v) \rvert \choose m}]\}$ for each $Y \in \P_e$, 
        \item $S_e = \{u_{e,M}: M$ is a monochromatic component of the unique $L|_{1_L}$-coloring of $G[1_L]$ intersecting $N_{G_e}^{\leq 1}[X_e]-X_e\}$, and 
        \item $S'_e = \{u'_{e,M}: M$ is a monochromatic component of the unique $L|_{1_L}$-coloring of $G[1_L]$ intersecting $N_{G_e}^{\leq 1}[X_e]-X_e\}$. 
    \end{itemize}

Let $T'$ be the tree obtained from $T_0$ by, for each $e \in U_E$, adding a new node $t_e$ adjacent to the end of $e$ in $T_0$.
For every $t \in V(T_0)$, let $X'_t = X_t$; for every $t \in V(T')-V(T_0)$, we know $t=t_e$ for some $e \in U_E$, and we let $X'_t = X_e \cup \bigcup_{Y \in \P_e}S_{e,Y} \cup S_e \cup S_e'$.
Let $\X'= (X'_t: t \in V(T'))$.
Clearly, $(T',\X')$ is a tree-decomposition of $H$ with adhesion at most $\max_{e \in U_E}\{\theta,\lvert X_e \rvert\} = \theta$.

If $\eta-1=0$, then let $(T'',\X'')=(T',\X' \cap (V(H)-Z))$.
If $\eta-1 \geq 1$, then let $t_0$ be a node of $T'$ with $X'_{t_0}-Z \neq \emptyset$ closest to $t^*$, let $v_0$ be a vertex in $X'_{t_0}-Z$, let $T''$ be the rooted tree obtained from $T'$ by adding a new node $t_0^*$ adjacent to $t^*$, where $t_0^*$ is the root of $T''$, and let $\X''=(X''_t: t \in V(T''))$, where $X''_{t_0^*}=\{v_0\}$, $X''_t=(X'_t-Z) \cup \{v_0\}$ if $t$ is in the path in $T'$ between $t^*$ and $t_0$, and $X''_{t}=X'_t-Z$ otherwise.

\begin{clm} \label{clm:construction}
$(T'',\X'')$ is an $(\eta-1,\theta)$-construction of $H-Z$.
\end{clm}

\begin{claimproof}
Note that $(T',\X' \cap (V(H)-Z))$ is a tree-decomposition of $H-Z$ with adhesion at most $\theta$.
Clearly, for every $tt' \in E(T')-E(T_0)$, one end of $tt'$ has no child, say $t'$, and every component of $H[(X'_{t'}-Z)-(X'_t-Z)] = H[\bigcup_{Y \in \P_e}S_{e,Y} \cup S_e \cup S_e']$ has at most two vertices.
By the definition of $T_0$, for every $tt' \in E(T_0)$, $X_t \cap X_{t'} \cap Z \neq \emptyset$, so $|(X'_t-Z) \cap (X_{t'}-Z)| \leq |X_t \cap X_t'|-1$.
Hence for every $tt' \in E(T_0)$, if $|(X'_t-Z) \cap (X_{t'}-Z)| \geq \eta-1$, then $|X_t \cap X_t'| \geq \eta$, so by (C1) for $(T,\X)$, one end of $tt'$, say $t'$, has no child in $T$ (and hence no child in $T'$), and every component of $H[(X'_{t'}-Z)-(X'_t-Z)] = G[X_{t'}-(X_t \cup Z)]$ has at most two vertices.
Therefore, $(T',\X' \cap V(H)-Z)$ satisfies (C1) for being an $(\eta-1,\theta)$-construction of $H-Z$.
Since $t^* \in V(T_0)$, $|X'_{t^*}-Z|=|X_{t^*}-Z| \leq \theta$.

If $\eta-1=0$, then $(T'',\X'')=(T',\X' \cap (V(H)-Z))$, so $(T'',\X'')$ is an $(\eta-1,\theta)$-construction of $H-Z$.
Hence we may assume $\eta-1 \geq 1$.
Then $(T'',\X'')$ is a tree-decomposition of $H-Z$ with adhesion at most $\theta$ such that $|X''_{t_0^*} \cap X''_{t^*}|=|\{v_0\}|=1 \leq \eta-1$, and for every $tt' \in E(T')$ with $X''_t \cap X''_{t'} \neq (X'_t - Z) \cap (X'_{t'}-Z)$, we have $|X''_t \cap X''_{t'}| \leq |\{v_0\}| = 1 \leq \eta-1$ (by the choice of $t_0$).
Since $(T',\X' \cap (V(H)-Z))$ satisfies (C1) for being an $(\eta-1,\theta)$-construction of $H-Z$, $(T'',\X'')$ satisfies (C1) for being an $(\eta-1,\theta)$-construction of $H$.
Moreover, $|X''_{t_0^*}|=|\{v_0\}|=1 \leq \eta-1 \leq \theta$.
So $(T'',\X'')$ is an $(\eta-1,\theta)$-construction of $H-Z$.
\end{claimproof}

Note that for each $e \in U_E$ and $Y \in \P_e$, there exists a bijection $\pi_{e,Y}$ from $\{v_{Y,i}: i \in [{\lvert \bigcup_{v \in V(G)}L(v) \rvert \choose m}]\}$ to the set of all $m$-element subsets of $\bigcup_{v \in V(G)}L(v)$.
We define a list-assignment $L_H$ of $H$ as follows:
    \begin{itemize}
        \item For each $e \in U_E$ and $Y \in \P_e$, define $L_H(v_{Y,i}) = \pi_{e,Y}(v_{Y,i})$ for each $i \in [{\lvert \bigcup_{v \in V(G)}L(v) \rvert \choose m}]$.
        \item For every $e \in U_E$ and $u_{e,M} \in S_e$, define $L_H(u_{e,M})$ to be the set consisting of the unique color on $M$.
        \item For every $e \in U_E$ and $u'_{e,M} \in S'_e$, define $L_H(u'_{e,M})=[m]$.
        \item For each $v \in V(H) \cap V(G)=V(G_0)$, define $L_H(v)=L(v)$.
    \end{itemize}

Our first goal is to prove that there exists an $L_H$-coloring of $H$ with smaller weak diameter (Claim \ref{claim4}), and then use it to obtain an $L$-coloring of $G$ with small weak diameter to finish the proof.
We will apply induction to $H$ to find a desired $L_H$-coloring, so we have to check the extendability of $L_H$, and the following a few claims are for this purpose.

\begin{clm} \label{clm:H_legitimate_1}
$L_H$ is list-assignment of $H$ such that the unique $L_H$-coloring of $H[1_{L_H}]$ has weak diameter in $H$ at most $k$.
\end{clm}

\begin{claimproof}
Let $W$ be a monochromatic component of the unique $L_H|_{1_{L_H}}$-coloring of $H[1_{L_H}]$.
To prove this claim, it suffices to show that $V(W)$ has weak diameter in $H$ at most $k$.

If $V(W) \cap V(G_0)=\emptyset$, then $W$ is contained in a component of $H-V(G_0)$, so $W$ has at most two vertices, and hence $V(W)$ has weak diameter in $H$ at most $1 \leq k$.
So we may assume $V(W) \cap V(G_0) \neq \emptyset$.
Note that $V(W) \cap V(G_0) \subseteq 1_{L_H} \cap V(G_0) \subseteq 1_L$.
Hence there exists a monochromatic component $M$ of the unique $L|_{1_L}$-coloring of $G[1_L]$ such that $V(W) \cap V(G_0) \cap V(M) \neq \emptyset$.
Let $x \in V(W) \cap V(G_0) \cap V(M)$.

Let $M_E = \{e \in U_E: V(M) \cap N_{G_e}^{\leq 1}[X_e]-X_e \neq \emptyset\}$.
Note that $u_{e,M}$ is defined for every $e \in M_E$. 
Let $O=\{u_{e,M}: e \in M_E\}$. 
We shall show that $V(W) \subseteq (V(M) \cap V(G_0)) \cup O$.

Suppose to the contrary that $V(W) - ((V(M) \cap V(G_0)) \cup O) \neq \emptyset$.
Since $x \in V(M) \cap V(G_0)$, there exists a path $P$ in $W$ from $x$ to a vertex $y \in V(W) - ((V(M) \cap V(G_0)) \cup O)$ such that $V(P)-\{y\} \subseteq (V(M) \cap V(G_0)) \cup O$.

We first suppose $y \not \in V(G_0)$.
Then by the definition of $L_H$, there exist $e \in U_E$ and $M''$ such that $y=u_{e,M''} \in V(P) \cap S_e$.
Since $y \not \in O$, we know $M'' \neq M$. 
Note that every vertex in $N_{G_e}^{\leq 1}[V(M'') \cap X_{T_e}-X_e]) \cap V(W)$ is adjacent in $G$ to $V(M'')$ and is in $V(G_0) \cap V(W) \subseteq 1_{L_H}$ with color equal to the color on $W$, which equals the color of $u_{e,M''}$ and the color of $M''$.
So $N_{G_e}^{\leq 1}[V(M'') \cap X_{T_e}-X_e]) \cap V(W) \subseteq V(M'')$.
Hence 
	\begin{align*}
		N_H(y) \cap V(P) \subseteq & N_H(u_{e,M''}) \cap V(G_0) \cap V(W) \\
		\subseteq & X_e \cap N_{G_e}^{\leq 1}[V(M'') \cap X_{T_e}-X_e] \cap V(W) \\
		\subseteq & X_e \cap V(M'') \\
		\subseteq & V(G_0) \cap V(M'').
	\end{align*}
Since $V(P)-\{y\} \subseteq (V(M) \cap V(G_0)) \cup O$, we have $M''=M$, a contradiction. 

Hence $y \in V(G_0)$.
Since $V(P)-\{y\} \subseteq (V(M) \cap V(G_0)) \cup O$, and for every $e \in M_E$, $N_H(u_{e,M}) \cap V(P) \subseteq V(M) \cap X_e$ by considering the color of $P \subseteq W$, we know that $y \in V(M) \cap V(G_0)$, a contradiction.

Therefore, $V(W) \subseteq (V(M) \cap V(G_0)) \cup O$.

For every $v \in V(W)$, if $v \in V(G_0)$, then let $\iota(v)=v$; otherwise, we know $v=u_{e_v,M}$ for some $e_v \in M_E$, and we let $\iota(v)$ be a vertex in $V(M) \cap X_{T_e}-X_e$.
Note that $\iota(v) \in V(M)$ for every $v \in V(W)$.

Let $a$ and $b$ be vertices in $V(W)$.
To prove this claim, it suffices to show that $\dist_H(a,b) \leq k$.
We may assume $a \neq b$, for otherwise we are done.

Since $L$ satisfies (L3), $V(M)$ has weak diameter in $G$ at most $k$, so there exists a path $R$ in $G$ between $\iota(a)$ and $\iota(b)$ with length at most $k$.
Since $a \neq b$ and $V(W) \subseteq (V(M) \cap V(G_0)) \cup O$, there exists no $e \in U_E$ such that both $\iota(a)$ and $\iota(b)$ are in $X_{T_e}-X_e$.
So for every $z \in \{a,b\}-V(G_0)$, we know that $\iota(z) \in V(M) \cap X_{T_{e_z}}-X_{e_z}$ and there exists a subpath $R_z$ of $R$ from $\iota(z)$ to $X_{e_z}$ internally disjoint from $X_{e_z}$; let $r_z$ be the end of $R_z$ in $X_{e_z}$.

For every $z \in \{a,b\}-V(G_0)$, we know $z=u_{e_z,M}$; if $|E(R_z)|=1$, then $r_z \in N_{G_{e_z}}^{\leq 1}[V(M) \cap X_{T_{e_z}}-X_{e_z}]$, so $u_{e_z,M}r_z \in E(H)$ and forms a path $\overline{R_z}$ in $H$ with length $1=|E(R_z)|$ between $z$ and $r_z$; if $|E(R_z)| \geq 2$, then since $|E(R_z)| \leq |E(R)| \leq k$, we know $\iota(z) \in (V(M) \cap X_{T_{e_z}}-X_{e_z}) \cap N_{G_{e_z}}^{\leq k}[\{r_z\}]$, so $r_z \in X_{e_z} \cap N_{G_{e_z}}^{\leq k}[V(M) \cap X_{T_{e_z}}-X_{e_z}] \subseteq N_H(u'_{e_z,M})$, and hence $u_{e_z,M}u'_{e_z,M}r_z$ is a path $\overline{R_z}$ in $H$ with length $2 \leq |E(R_z)|$ between $z$ and $r_z$. 

For every subpath $R'$ of $R$ between two vertices in $X_e$ (for some $e \in U_E$) contained in $G_e$ internally disjoint from $X_e$ with $V(R')-X_e \neq \emptyset$, the ends of $R'$ are contained in the same part $Y$ of $\P_e$ (since $|E(R')| \leq |E(R)| \leq k$), so there exists a path $\overline{R'}$ in $H$ between the ends of $R'$ with exactly one internal vertex $v_{Y,1}$; note that $|E(\overline{R'})| \leq |E(R')|$.

Let $\overline{R}$ be the walk obtained from $R$ by, for each $z \in \{a,b\}-V(G_0)$, replacing $R_z$ by $\overline{R_z}$, and by, for each subpath $R'$ of $R$ between two vertices in $X_e$ (for some $e \in U_E$) contained in $G_e$ internally disjoint from $X_e$ with $V(R')-X_e \neq \emptyset$, replacing $R'$ by $\overline{R'}$.
Then $\overline{R}$ is a walk in $H$ between $a$ and $b$ such that $|E(\overline{R})| \leq |E(R)| \leq k$.
So $\dist_H(a,b) \leq k$.
\end{claimproof}

For every $v \in V(H)-Z$,
    \begin{itemize}
        \item if $v \in N_H^{\leq k+r}[Z]$, then let $L_H'(v)$ be a superset of $L_H(v)$ with size $m$,
        \item otherwise, let $L_H'(v)=L_H(v)$.
    \end{itemize}

\begin{clm} \label{clm:H_legitimate_2}
$L_H'$ is a $(T'',\X'',m,s,r,k)$-legitimate list-assignment of $H-Z$.
\end{clm}

\begin{claimproof}
Clearly, $|L_H'(v)| \in \{1,m\}$ for every $v \in V(H)-Z$, so (L1) holds.

Now we show that $L_H'$ satisfies (L2).
Let $t \in V(T'')$.
If $|X''_t| \leq 1$, then $1_{L_H'} \cap X''_t$ is $(1,0)$-centered (and hence $(s,r)$-centered) in $(H-Z)[X''_t]$.
So we may assume $|X''_t| \geq 2$.
In particular, $t \in V(T')$.
If $t \in V(T')-V(T_0)$, then $t=t_e$ for some $e \in U_E$, so $1_{L_H'} \cap X''_{t_e} \subseteq X''_{t_e} \subseteq N_{(H-Z)[X''_{t_e}]}^{\leq 1}[X_e]$ is $(\theta,1)$-centered (and hence $(s,r)$-centered since $s \geq \theta$) in $(H-Z)[X''_{t_e}]$.
So we may assume $t \in V(T_0)$ with $|X''_t| \geq 2$.
Hence $1_{L_H'} \cap X''_t \subseteq 1_L \cap X_t-Z$.
Since $L$ is hereditarily $(T,\X,m,s,r,k,N)$-extendable, $L$ is $(T,\X,m,s,r,k)$-legitimate, so $1_{L_H'} \cap X''_t \subseteq 1_L \cap X_t-Z \subseteq 1_L \cap X_t$ is $(s,r)$-centered in $G[X_t]$.
So there exists $S \subseteq X_t$ with $|S| \leq s$ such that $1_{L_H'} \cap X''_t \subseteq N_{G[X_t]}^{\leq r}[S]$.
To show that $L_H'$ satisfies (L2), it suffices to show $1_{L_H'} \cap X''_t \subseteq N_{(H-Z)[X''_t]}^{\leq r}[S-Z]$.
Let $x \in 1_{L_H'} \cap X''_t$.
Since $1_{L_H'} \cap X''_t \subseteq N_{G[X_t]}^{\leq r}[S]$, there exists a path $P_x$ in $G[X_t]$ between $x$ and $S$ with $|E(P_x)| \leq r$.
Note that $(H-Z)[X''_t] = G[X_t]-Z$.
If $V(P_x) \cap Z=\emptyset$, then $P_x \subseteq G[X_t]-Z = (H-Z)[X''_t]$, so $x \in N_{(H-Z)[X''_t]}^{\leq r}[S-Z]$.
If $V(P_x) \cap Z \neq \emptyset$, then $x \in N_{G[X_t]}^{\leq r}[Z] \subseteq N_H^{\leq r}[Z]$ (since $G[X_t]=H[X_t] \subseteq H$), so $x \not \in 1_{L_H'}$ by the definition of $L_H'$, a contradiction.
This shows that $L_H'$ satisfies (L2).

Finally, we show that $L_H'$ satisfies (L3).
Let $W$ be a monochromatic component of the unique $L_H'|_{1_{L_H'}}$-coloring of $(H-Z)[1_{L_H'}]$.
To prove this claim, it suffices to show that $V(W)$ has weak diameter in $H-Z$ at most $k$.
Let $a,b \in V(W)$.
It suffices to show $\dist_{H-Z}(a,b) \leq k$.

Since $L_H'(v) \supseteq L_H(v)$ for every $v \in V(H)-Z$, we know $1_{L_H'} \subseteq 1_{L_H}$.
So $V(W)$ is contained in a monochromatic component $M$ of the unique $L_H|_{1_{L_H}}$-coloring of $H[1_{L_H}]$.
By Claim \ref{clm:H_legitimate_1}, $V(M)$ has weak diameter in $H$ at most $k$.
So there exists a path $P$ in $H$ between $a$ and $b$ with $|E(P)| \leq k$.
If $V(P) \subseteq V(H)-Z$, then $\dist_{H-Z}(a,b) \leq |E(P)| \leq k$ and we are done.
So we may assume $V(P) \cap Z \neq \emptyset$.
Hence $a \in N_H^{\leq |E(P)|}[Z] \subseteq N_H^{\leq k}[Z]$.
By the definition of $L_H'$, we know $|L_H'(a)|=m$, so $a \not \in 1_{L_H'} \supseteq V(W)$, a contradiction.
\end{claimproof}

\begin{clm} \label{clm:H_extendable}
$L_H'$ is a hereditarily $(T'',\X'',m,s,r,k,N_1+f_2(N))$-extendable list-assignment of $H-Z$.
\end{clm}

\begin{claimproof}
By Claim \ref{clm:H_legitimate_2}, it suffices to verify (H2) for $L_H'$.
Let $W$ be an induced subgraph of $H-Z$.
Let $L_W$ be a $(T'',\X'' \cap V(W),m,s,r,k)$-legitimate list-assignment of $W$ with $L_W(v) \supseteq L_H'(v)$ for every $v \in V(W)$.
It suffices to verify (E2) for $L_W$ being $(T'',\X'' \cap V(W), m,s,r,k,N_1+f_2(N))$-extendable.

Let $t \in V(T'')$.
Let $G_t$ be a child-extension at $t$ with respect to $(T'',\X'' \cap V(W))$.
Let $L_t$ be a list-assignment of $G_t$ with $L_t|_{X''_t \cap V(W)} = L_W|_{X''_t \cap V(W)}$.
It suffices to show that there exists an $L_t$-coloring of $G_t$ with weak diameter in $G_t$ at most $N_1+f_2(N)$.

If $|X''_t| \leq 1$, then $|X_t'' \cap V(W)| \leq 1$, implying that every component of $G_t$ has diameter at most 4, so there exists an $L_t$-coloring of $G_t$ with weak diameter in $G_t$ at most $4 \leq N \leq f_2(N)$.
Hence we may assume $|X''_t| \geq 2$.
In particular, $t \in V(T')$.
We first assume $t \in V(T')-V(T_0)$.
Then $t=t_e$ for some $e \in U_E$.
So $t$ has no child.
Hence $G_t = W[X''_t \cap V(W)]$, and for every component $C$ of $G_t$, if $V(C) \cap X_e =\emptyset$, then $C$ has at most two vertices; otherwise $V(C) \subseteq N_{W[X''_t \cap V(W)]}^{\leq 2}[X_e \cap V(W)]$.
So the vertex-set of each component of $G_t$ is $(\theta,2)$-centered in $G_t$.
Hence there exists an $L_t$-coloring of $G_t$ with weak diameter in $G_t$ at most $N_1$ by Lemma \ref{all_centered}. 

So we may assume $t \in V(T_0)$ with $|X''_t| \geq 2$.
Hence $X''_t=X_t$ by the choice of $t_0$.
Since $t \in V(T_0)$, we have $X''_t \cap V(W) \subseteq V(G_0)$, so for every $v \in X''_t \cap V(W)$, $L_W(v) \supseteq L_H'(v) \supseteq L_H(v) \supseteq L(v)$.
Let $G_W=G[V(W) \cap X''_t]$.
Let $L_{G_W}$ be a list-assignment of $G_W$ such that $L_{G_W}(v)$ is a superset of $L_W(v)$ with size $m$ for every $v \in X''_t \cap V(W)$.
So $L_{G_W}(v) \supseteq L(v)$ for every $v \in V(G_W)$.
Since $1_{L_{G_W}}=\emptyset$, $L_{G_W}$ is $(T,\X \cap V(G_W),m,s,r,k)$-legitimate.
Since $L$ is hereditarily $(T,\X,m,s,r,k,N)$-extendable, by (H2), $L_{G_W}$ is $(T,\X \cap V(G_W),m,s,r,k,N)$-extendable.
Note that the bag at $t$ in $(T'',\X'' \cap V(W))$ equals the bag at $t$ in $(T,\X \cap V(G_W))$.
And for every child $t'$ of $t$ in $T''$, either $t'$ is a child of $t$ in $T$, or $t'=t_e$ for some $e \in U_E$.
For the former, the intersection of the bags of $t$ and $t'$ in $(T'',\X'' \cap V(W))$ equals the intersection of the bags of $t$ and $t'$ in $(T,\X \cap V(G_W))$; for the latter, the intersection of the bags of $t$ and $t'$ in $(T'',\X'' \cap V(W))$ equals $X_e \cap V(G_W)$.
So $G_t$ is also a child-extension at $t$ with respect to $(T,\X \cap V(G_W))$.
Let $L_t'$ be the list-assignment of $G_t$ such that $L_t'|_{X''_t \cap V(W)} = L_{G_W}|_{X''_t \cap V(W)}$ and $L_t'|_{V(G_t)-(X''_t \cap V(W))} = L_t|_{V(G_t)-(X''_t \cap V(W))}$.
Since $L_{G_W}$ is $(T,\X \cap V(G_W),m,s,r,k,N)$-extendable, there exists a $L_t'$-coloring $c_t$ of $G_t$ with weak diameter in $G_t$ at most $N$.

Note that for every $v \in V(G_t)$, if $L_t'(v) \neq L_t(v)$, then $v \in 1_{L_t} \cap X''_t \cap V(W) = 1_{L_W} \cap X''_t \cap V(W)$.
Since $L_W$ is $(T'',\X'' \cap V(W),m,s,r,k)$-legitimate, by (L2), the set $\{v \in V(G_t): L_t'(v) \neq L_t(v)\} \subseteq 1_{L_W} \cap X''_t \cap V(W)$ is $(s,r)$-centered in $W[X_t'' \cap V(W)]$ and hence $(s,r)$-centered in $G_t$.
Let $Z'=\{v \in V(G_t): L_t'(v) \neq L_t(v)\}$.
Let $c_{Z'}$ be the unique $L_t|_{Z'}$-coloring of $G_t[Z']$.
Applying Lemma \ref{add_centered} by taking $(k,r,N,G,Z,m,R,c_Z,c)=(s,r,N,G_t,Z',|\bigcup_{v \in V(G_t)}L_t(v)|,\emptyset,c_{Z'},c_t|_{V(G_t)-Z'})$, the coloring $c_t|_{V(G_t)-Z'} \cup c_{Z'}$ is an $L_t$-coloring of $G_t$ with weak diameter in $G_t$ at most $f_2(N)$.
This proves the claim.
\end{claimproof}

\begin{clm} \label{claim4}
There is an $L_H$-coloring $c_H$ of $H$ with weak diameter in $H$ at most $f_3(f_1(f^*(\eta-1,N_1+f_2(N))))$ such that $c_H(v)=c_Z(v)$ for every $v \in Z$. 
\end{clm}

\begin{claimproof}
By Claim \ref{clm:construction}, $H-Z$ is $(\eta-1,\theta)$-constructible with an $(\eta-1,\theta)$-construction $(T'',\X'')$.
By Claim \ref{clm:H_extendable}, $L_H'$ is a hereditarily $(T'',\X'',m,s,r,k,N_1+f_2(N))$-extendable list-assignment of $H-Z$.
By the minimality of $\eta$, we know that the unique $L_H'|_\emptyset$-coloring of the empty graph can be extended to an $L_H'$-coloring $c_H'$ of $H-Z$ with weak diameter in $H-Z$ at most $f^*(\eta-1,N_1+f_2(N))$.
Since $H-Z \subseteq H$, $c_H'$ has weak diameter in $H$ at most $f^*(\eta-1,N_1+f_2(N))$.

By Claim \ref{claim_Zbasic}, $Z=N_G^{\leq (k+2)(\theta+1)}[X_{t^*}]$.
Since $Z \subseteq V(H)$, $Z=N_H^{\leq (k+2)(\theta+1)}[X_{t^*}]$ is $(\theta,(k+2)(\theta+1))$-centered in $H$.
By Lemma \ref{add_centered}, $c_H' \cup c_Z$ is a coloring of $H$ with weak diameter in $H$ at most $f_1(f^*(\eta-1,N_1+f_2(N)))$ such that $(c_H' \cup c_Z)|_Z=c_Z$.

Let $Z'=\{v \in V(H)-Z: L_H'(v) \neq L_H(v)\}$.
Note that $Z' \subseteq N_H^{\leq k+r}[Z]-Z \subseteq N_H^{\leq (k+2)(\theta+1)+k+r}[X_{t^*}]-Z$ is $(\theta,(k+2)(\theta+1)+k+r)$-centered in $H$.
Let $c_{Z'}$ be an $L_H|_{Z'}$-coloring of $H[Z']$.
Then by Lemma \ref{add_centered}, $(c_H' \cup c_Z)|_{V(H)-Z'} \cup c_{Z'}$ is a coloring of $H$ with weak diameter in $H$ at most $f_3(f_1(f^*(\eta-1,N_1+f_2(N))))$.
Note that $(c_H' \cup c_Z)|_{V(H)-Z'} \cup c_{Z'}$ is an $L_H$-coloring such that $((c_H' \cup c_Z)|_{V(H)-Z'} \cup c_{Z'})|_Z = (c_H' \cup c_Z)|_Z=c_Z$.
This proves the claim.
\end{claimproof}
	
For every $e \in U_E$, 
    \begin{itemize}
        \item let $\sigma_e$ be an ordering of the vertices in $X_e$; that is, $\sigma_e$ is a bijection from $X_e$ to $[\lvert X_e \rvert]$;
        \item for each $i\in [\theta+2]$, define $Z_{i,e} = N_{G_e}^{\leq (k+2)i}[X_e]$.
    \end{itemize}
Let $Z_0 = \bigcup_{e \in U_E} X_e$.
Recall that $G_1=\bigcup_{e\in U_E} G_e$.
So $Z_0=V(G_0)\cap V(G_1)$.
For $i \in [\theta+2]$, let $Z_i = \bigcup_{e\in U_e} Z_{i,e}$.
	
\begin{clm} \label{claim5}
For any $e \in U_E$, $i \in [\theta+1]$ and $v \in Z_i \cap X_{T_e}-X_e$, there exists $Y \in \P_e$ such that $v \in N_{G_e}^{\leq (k+2)i}[Y]$, and for every $Y' \in \P_e-\{Y\}$, $v \not \in N_{G_e}^{\leq (k+2)i}[Y']$.
\end{clm}

\begin{claimproof}
Since $v \in Z_{i} \cap X_{T_e}-X_e$, there exists a path $P$ in $G_1$ from $v$ to $X_e$ internally disjoint from $X_e$ of length at most $(k+2)i$.
Since $P$ is internally disjoint from $X_e$, $P$ is a path in $G_e$.
Let $y$ be the vertex in $V(P) \cap X_e$.
Let $Y$ be the member of $\P_e$ containing $y$.
So $v \in N_{G_e}^{\leq (k+2)i}[Y]$.
Let $Y'$ be any member of $\P_e-\{Y\}$.
If $v \in N_{G_e}^{\leq (k+2)i}[Y']$, then there exists a walk in $G_e$ from $Y$ to $Y'$ of length at most $2(k+2)i \leq 2(k+2)(\theta+1)+1$, so $Y=Y'$ by the definition of $\P_e$, a contradiction.
Hence $v \not \in N_{G_e}^{\leq (k+2)i}[Y']$.
\end{claimproof}
	
For every $v \in Z_1-V(H) = Z_1-Z_0 \subseteq Z_{\theta+1}-Z_0$, there is a unique $e_v\in U_E$ with $v \in Z_{\theta+1} \cap X_{T_{e_v}}-X_{e_v}$, hence by Claim \ref{claim5}, there exists a unique pair $(e_v,Y_v)$ with $e_v \in U_E$ and $Y_v \in \P_{e_v}$ such that $v \in N_{G_{e_v}}^{\leq (k+2)(\theta+1)}[Y_v]$; moreover, if $v \not \in 1_L$, then $|L(v)|=m$ by (L1), so there exists $s_v \in S_{e_v,Y_v}$ such that $L(v)=\pi_{e_v,Y_v}(s_v)=L_H(s_v)$ by the definition of $L_H$.
	
Let $c_H$ be the $L_H$-coloring of $H$ mentioned in Claim \ref{claim4}. 
This assigns a color to all vertices in $V(H)$ and in particular all vertices in $V(G_0)$. 
Our next goal is to extend this coloring of $G_0$ to an $L$-coloring of $G$ and then show that this coloring has small weak diameter to finish the proof.

Let $c_{\theta+1}$ be a function with domain $Z_{\theta+1}$ such that 
	\begin{itemize}
		\item $c_{\theta+1}(u)=c_H(u)$ for every $u \in Z_1 \cap V(H) = Z_0$, 
        \item $c_{\theta+1}(u)$ is the unique element in $L(u)$ for every $u \in Z_{\theta+1} \cap 1_L$,
		\item $c_{\theta+1}(u) = c_H(s_u)$ for every $u \in Z_1-(V(H) \cup 1_L) = Z_1-(Z_0 \cup 1_L)$, and
		\item for every $u \in Z_{\theta+1}-(Z_1 \cup 1_L)$, if $i_u$ is the element in $[\theta]$ such that $u \in Z_{i_u+1}-Z_{i_u}$, then
		\begin{itemize}
			\item when $\lvert X_{e_u} \rvert \geq i_u$, $c_{\theta+1}(u)$ is an element in $L(u)-\{c_H(\sigma_{e_u}^{-1}(i_u))\}$, and
			\item when $\lvert X_{e_u} \rvert < i_u$, $c_{\theta+1}(u)$ is an element in $L(u)$.
		\end{itemize}
	\end{itemize}
Note that $c_{\theta+1}$ is well-defined, since $L_H|_{Z_0}=L|_{Z_0}$ and $m \geq 2$.
Moreover, $c_{\theta+1}$ is an $L|_{Z_{\theta+1}}$-coloring of $G[Z_{\theta+1}]$.
	
For every $e \in U_E$, let $Z_e = Z_{\theta+1} \cap V(G_e)$, and let $c_e$ be the function with domain $Z_e$ such that $c_e(v)=c_{\theta+1}(v)$ for every $v \in Z_e$.
By Claim \ref{claim5}, for every $e \in U_E$, $Z_e \subseteq N_{G_e}^{\leq (k+2)(\theta+1)}[X_e]$.
	
\begin{clm} \label{claim 6}
For every $e \in U_E$, $c_e$ can be extended to an $L|_{V(G_e)}$-coloring $c_e'$ of $G_e$ with weak diameter in $G_e$ at most $f^*(\eta,N)$.
\end{clm}
	
\begin{claimproof}
Let $e \in U_E$.
Let $T_e'$ be the rooted tree obtained from $T_e$ by adding a new node $t_e^*$ adjacent to the end of $e$ in $V(T_e)$, where $t_e^*$ is the root of $T_e'$.
Let $X^e_{t_e^*}=X_e$.
For every $t \in V(T_e)$, let $X^e_t=X_t$.
Let $\X^e=(X^e_t: t \in V(T_e'))$.
Then $(T'_e,\X^e)$ is a rooted tree-decomposition of $G_e$ with adhesion at most $\theta$ satisfying (C1) for being an $(\eta,\theta)$-construction and $|X^e_{t_e^*}|=|X_e| \leq \theta$.
By Claim \ref{clm:X_e_nonempty}, $X_e \neq \emptyset$.
So $X^e_{t_e^*} \neq \emptyset$.
Hence $(T_e',\X^e)$ is an $(\eta,\theta)$-construction of $G_e$.

Let $L_e$ be the list-assignment of $G_e$ such that $L_e|_{V(G_e)-Z_e}=L|_{V(G_e)-Z_e}$, and for every $v \in Z_e$, $L_e(v)$ is a superset of $L(v)$ with size $m$.

We first show that $L_e$ is $(T_e',\X^e,m,s,r,k)$-legitimate.
Clearly, (L1) holds for $L_e$.
Since $L$ is $(T,\X,m,s,r,k)$-legitimate, for every $t \in V(T_e)$, $1_{L_e} \cap X^e_t=1_{L_e} \cap X_t \subseteq 1_L \cap X_t$ is $(s,r)$-centered in $G[X_t]=G_e[X_t]$ by (L2) for $L$.
And $1_{L_e} \cap X^e_{t^*_e} = 1_{L_e} \cap X_e \subseteq 1_{L_e} \cap Z_e = \emptyset$.
So (L2) holds for $L_e$.
Let $W$ be a monochromatic component of the unique $L_e|_{1_{L_e}}$-coloring of $G_e[1_{L_e}]$.
To show that (L3) holds for $L_e$, it suffices to show that $V(W)$ has weak diameter in $G_e$ at most $k$.
Since $1_{L_e} \subseteq 1_L$, $W$ is contained in a monochromatic component $M$ of the unique $L|_{1_L}$-coloring of $G[1_L]$.
Since $L$ is $(T,\X,m,s,r,k)$-legitimate, $V(M)$ has weak diameter in $G$ at most $k$.
Let $a,b$ be vertices in $W$.
Since $V(W) \subseteq V(M)$, there exists a path $P$ in $G$ between $a$ and $b$ with length at most $k$.
Since $Z_e \cap 1_{L_e}=\emptyset$, $V(W) \cap Z_e=\emptyset$.
So $\{a,b\} \cap Z_e=\emptyset$.
Hence $\{a,b\} \subseteq X_{T_e}-Z_e$.
If $P \not \subseteq G_e$, then $V(P) \cap X_e \neq \emptyset$, so $a \in N_{G_e}^{\leq |E(P)|}[X_e] \subseteq N_{G_e}^{\leq k}[X_e] \subseteq Z_1 \cap V(G_e) \subseteq Z_e$, a contradiction.
So $P \subseteq G_e$, and hence $\dist_{G_e}(a,b) \leq k$.
This shows that $W$ has weak diameter in $G_e$ at most $k$.

Therefore, $L_e$ is $(T_e',\X^e,m,s,r,k)$-legitimate. 
Note that $L_e(v) \supseteq L(v)$ for every $v \in V(G_e)$; for every $t \in V(T)-V(T_e)$, we know $1_{L_e} \cap X_t \subseteq 1_{L_e} \cap V(G_e) \cap X_t \subseteq 1_{L_e} \cap X_e \subseteq 1_{L_e} \cap Z_e = \emptyset$.
Since $L_e$ is $(T_e',\X^e,m,s,r,k)$-legitimate, by Lemma \ref{hereditary_list_2}, $L_e$ is hereditarily $(T_e',\X^e,m,s,r,k,N)$-extendable.

Recall that $(T_e',\X^e)$ is an $(\eta,\theta)$-construction of $G_e$, and $L_e$ is hereditarily $(T_e',\X^e,m,s,r,k, \allowbreak N)$-extendable.
Note that $Z_e \subseteq N_{G_e}^{\leq (k+2)(\theta+1)}[X_e] = N_{G_e}^{\leq (k+2)(\theta+1)}[X^e_{t^*_e}]$ and $c_e$ is an $L_e|_{Z_e}$-coloring of $G_e[Z_e]$.
Since $(\eta,|V(G_e)-Z_e|+|V(G_e)|)$ is lexicographically smaller than $(\eta,|V(G)-Z|+|V(G)|)$, $c_e$ can be extended to an $L_e$-coloring $c_e'$ of $G_e$ with weak diameter in $G_e$ at most $f^*(\eta,N)$.
Since $\{v \in V(G_e): L_e(v) \neq L(v)\} \subseteq Z_e$ and $c_e'|_{Z_e}=c_e$, we know that $c_e$ is an $L|_{V(G_e)}$-coloring of $G_e$.
\end{claimproof}
	
For every $e \in U_E$, let $c_e'$ be the $L|_{V(G_e)}$-coloring mentioned in Claim \ref{claim 6}.
Define $$c = c_H|_{V(G_0)} \cup \bigcup_{e \in U_E}c'_e.$$
If $v\in V(G_0)\cap V(G_e)$ for some $e\in U_E$, then $c_e'(v)=c_{\theta+1}(v)=c_H(v)$ and hence  $c$ is well-defined. 
Moreover, $c|_Z=c_H|_Z=c_Z$ and $c_H(v) \in L(v)$ for every $v \in V(G_0)$.
So $c$ is an $L$-coloring of $G$ that can be obtained by extending $c_Z$. 
	
Therefore, there exists a $c$-monochromatic component $M$ in $G$ such that the weak diameter in $G$ of $V(M)$ is greater than $f^*(\eta,N)$.

\begin{clm} \label{claim_intersect}
For every $e\in U_E$, if $V(M) \cap V(G_e)-X_e \neq \emptyset$, then $V(M) \cap X_e \neq \emptyset$. 
\end{clm} 
	
\begin{claimproof}
Suppose to the contrary that there exists $e \in U_E$ such that $V(M) \cap V(G_e)-X_e \neq \emptyset$ and $V(M) \cap X_e=\emptyset$.
Then $V(M)\subseteq V(G_e)-X_e$.
Hence $M \subseteq G_e$. 
So $M$ is a $c_e'$-monochromatic component in $G_e$.
Hence $V(M)$ has weak diameter in $G_e$ (and hence in $G$) at most $f^*(\eta,N)$.
\end{claimproof}
	
\begin{clm} \label{claim_oneZ}
$V(M)$ does not intersect both $V(G_0)$ and $Z_{\theta+2}- Z_{\theta+1}$. 
\end{clm} 
	
\begin{claimproof}
Suppose that $V(M)$ intersects both $V(G_0)$ and $Z_{\theta+2}-Z_{\theta+1}$.
Since $M$ is a connected subgraph of $G$, there exists $e \in U_E$ such that $V(M)$ intersects $X_e$ and $V(G_e)-X_e$, and there exists a path $P$ in $M$ from a vertex $v_M$ in $X_e$ intersecting $(Z_{i+1}- Z_i) \cap X_{T_e}$ for every $i \in [\theta+1]$ internally disjoint from $X_e$. 
So for every $i \in [\theta+1]$, there exists a subpath $P_i$ of $P$ contained in $P[V(P) \cap Z_{i+1}-Z_i]$ intersecting both $N_{G_e}^{\leq (k+2)i+1}[X_e]-N_{G_e}^{\leq (k+2)i}[X_e]$ and $N_{G_e}^{\leq (k+2)(i+1)}[X_e]-N_{G_e}^{\leq (k+2)i+k+1}[X_e]$.
Since $|X_e| \leq \theta$, we know that $P_{\sigma_e(v_M)}$ is defined, and we denote it by $P'$.
For every $v \in V(P') \subseteq (Z_{\sigma_e(v_M)+1}-Z_{\sigma_e(v_M)}) \cap X_e$, if $v \not \in 1_L$, then by the definition of $c$, $c(v)=c_{\theta+1}(v) \neq c_H(v_M)=c(v_M)$, contradicting that $M$ is a $c$-monochromatic component.
So $V(P') \subseteq 1_L$.
Hence $P'$ is contained in a monochromatic component $M_{P'}$ of the unique $L|_{1_L}$-coloring of $G[1_L]$ and contains a vertex $a \in N_{G_e}^{\leq (k+2)\sigma_e(v_M)+1}[X_e]-N_{G_e}^{\leq (k+2)\sigma_e(v_M)}[X_e]$ and a vertex $b \in N_{G_e}^{\leq (k+2)(\sigma_e(v_M)+1)}[X_e]-N_{G_e}^{\leq (k+2)\sigma_e(v_M)+k+1}[X_e]$.
Since $L$ is $(T,\X,m,s,r,k)$-legitimate, $M_{P'}$ has weak diameter in $G$ at most $k$ by (L3), so there exists a path $Q$ in $G$ with $|E(Q)| \leq k$ between $a$ and $b$.
Since $\sigma_e(v_M) \geq 1$, we have $Q \subseteq G_e$, so $\dist_{G_e}(a,b) \leq k$.
Hence $b \in N_{G_e}^{\leq \dist_{G_e}(a,b)}[\{a\}] \subseteq N_{G_e}^{\leq (k+2)\sigma_e(v_M)+1+k}[X_e]$, a contradiction.
\end{claimproof}

Let $$B = \{e \in U_E: V(M) \cap V(G_e)-X_e \neq \emptyset\}.$$
	
\begin{clm} \label{claim_G0_and_others}
$V(M) \cap V(G_0) \neq \emptyset$ and $V(M) \subseteq V(G_0) \cup \bigcup_{e \in B}(X_{T_e} \cap Z_{\theta+1})$.
\end{clm}
	
\begin{claimproof}
If $V(M) \cap V(G_0)=\emptyset$, then there exists $e^* \in U_E$ such that $V(M) \cap V(G_{e^*})-X_{e^*} \neq \emptyset = V(M) \cap V(G_0) \supseteq V(M) \cap X_{e^*}$, contradicting Claim \ref{claim_intersect}.
So $V(M) \cap V(G_0) \neq \emptyset$.

To prove this claim, it suffices to show $V(M)-V(G_0) \subseteq \bigcup_{e \in B}(X_{T_e} \cap Z_{\theta+1})$.
Let $v \in V(M)-V(G_0)$.
So there exists $e \in U_E$ such that $v \in V(M) \cap V(G_e)-V(G_0)=V(M) \cap V(G_e)-X_e$.
By the definition of $B$, $e \in B$.
Since $V(M) \cap V(G_0) \neq \emptyset$, by Claim~\ref{claim_oneZ}, $V(M) \cap Z_{\theta+2}-Z_{\theta+1} = \emptyset$.
Since $M$ is connected and $V(M) \cap V(G_0) \neq \emptyset$, $v \in X_{T_e} \cap Z_{\theta+1}$.
This proves the claim.
\end{claimproof}

For every $x \in V(M)-V(G_0)$, by Claim \ref{claim_G0_and_others}, there exists a unique $e \in B$ with $x \in X_{T_e} \cap Z_{\theta+1}-X_e$; and recall that there exists a unique $Y_x \in \P_e$ such that $x \in N_{G_e}^{\leq (k+2)(\theta+1)}[Y_x]$ by Claim \ref{claim5}.

\begin{clm} \label{clm:path_same_Y}
For every $e \in U_E$ and every path $P=x_1x_2...x_{|V(P)|}$ in $M \cap G_e-X_e$, we have $Y_{x_1}=Y_{x_2}=...=Y_{x_{|V(P)|}}$.
\end{clm}

\begin{claimproof}
By Claim \ref{claim_G0_and_others}, $V(P) \subseteq Z_{\theta+1} \cap X_{T_e}-X_e$.
For each $v \in V(P)$, since $v \in N_{G_e}^{\leq (k+2)(\theta+1)}[Y_{v}]$, there exists a path $Q_v$ in $G_e$ with $|E(Q_v)| \leq (k+2)(\theta+1)$ between $Y_v$ and $v$.
Hence for any $i \in [\lvert V(P) \rvert-1]$, if $Y_{x_i} \neq Y_{x_{i+1}}$, then $(Q_{x_i} \cup Q_{x_{i+1}})+x_ix_{i+1}$ is a walk in $G_e$ from $Y_{x_i}$ to $Y_{x_{i+1}}$ with length at most $2(k+2)(\theta+1)+1$, so $Y_{x_i}=Y_{x_{i+1}}$ by the definition of $\P_e$, a contradiction.
So $Y_{x_i}=Y_{x_{i+1}}$ for every $i \in [\lvert V(P) \rvert-1]$.		
\end{claimproof}

\begin{clm} \label{clm:path_Gbranch_H}
For every $e \in U_E$ and every path $P$ in $M \cap G_e$ with both ends in $X_e$ internally disjoint from $X_e$, there exists a $c_H$-monochromatic path $\tilde{P}$ in $H$ having the same ends as $P$.
\end{clm}

\begin{claimproof}
Let $a$ and $b$ be the ends of $P$.
Since $\{a,b\} \subseteq X_e \subseteq V(H) \cap V(G_0)$, we have $c_H(a)=c(a)=c(M)=c(b)=c_H(b)$.
We may assume $a \neq b$ and $P$ contains at least one internal vertex, for otherwise we are done by taking $\tilde{P}=P$.
Let $Q=P-\{a,b\}$.

We first assume $V(Q) \subseteq 1_L$.
Then $Q$ is contained in a monochromatic component $M_Q$ of the unique $L|_{1_L}$-coloring of $G[1_L]$.
Since $\{a,b\} \subseteq X_e$ and $P$ is internally disjoint from $X_e$, the ends of $Q$ are contained in $N_{G_e}^{\leq 1}[X_e]-X_e$.
So $u_{e,M_Q}$ is defined.
Note that $\{a,b\} \subseteq X_e \cap N_{G_e}^{\leq 1}[V(Q)] \subseteq X_e \cap N_{G_e}^{\leq 1}[V(M_Q) \cap X_{T_e}-X_e] \subseteq N_H(u_{e,M_Q})$.
By the definition of $L_H$, we have $c_H(u_{e,M_Q})=c(M)$.
So $au_{e,M_Q}b$ is a desired $c_H$-monochromatic path $\tilde{P}$ in $H$.

So we may assume $V(Q)-1_L \neq \emptyset$.
Let $q$ be the vertex in $V(Q)-1_L$ such that the subpath $P_q$ of $P$ between $a$ and $q$ is as short as possible.

Suppose to the contrary that $q \not \in Z_1$.
Let $R=P_q-\{a,q\}$.
Since $q \not \in Z_1$, $V(R) \neq \emptyset$.
By the minimality of $|V(P_q)|$, $V(R) \subseteq 1_L$.
So $R$ is contained in a monochromatic component $M_R$ of the unique $L|_{1_L}$-coloring of $G[1_L]$.
Since $q \in X_{T_e}-Z_1$, $q \in Z_{\theta+1} \cap X_{T_e}-Z_1$ by Claim \ref{claim_G0_and_others}.
So $R$ contains a vertex $r_1 \in N_{G_e}^{\leq 1}[X_e]-X_e$ and a vertex $r_2 \in N_{G_e}^{\leq k+2}[X_e]-N_{G_e}^{\leq k+1}[X_e]$.
Then $\dist_G(r_1,r_2) \geq k+1$.
So the weak diameter in $G$ of $V(M_R)$ is greater than $k$, contradicting that $L$ is $(T,\X,m,s,r,k)$-legitimate.

Hence $q \in Z_1 \cap X_{T_e}-X_e$.
Recall that $s_q$ is a vertex in $S_{e,Y_q}$ such that $L(q)=L_H(s_q)$.
Since $q \in Z_1-(X_e \cup 1_L)$, $c_{\theta+1}(q)=c_H(s_q)$.
Hence by the definition of $c$, $c(M)=c(q)=c_{\theta+1}(q)=c_H(s_q)$.

Denote $Q$ by $x_1x_2\dots x_{\lvert V(Q) \rvert}$, where $ax_1 \in E(P)$ and $bx_{|V(Q)|} \in E(P)$. 
By Claim \ref{clm:path_same_Y}, $Y_{x_1}=Y_{x_{|V(Q)|}}=Y_q$.
By the uniqueness of $Y_{x_1}$ stated in Claim 10, $a \in Y_{x_1}$; similarly, $b \in Y_{x_{|V(Q)|}}$.
Since $Y_{x_1}=Y_{x_{|V(Q)|}}=Y_q$, we know that $s_q \in S_{e,Y_q}$ is adjacent in $H$ to both $a$ and $b$.
Hence $as_qb$ is a desired $c_H$-monochromatic path $\tilde{P}$ in $H$.
\end{claimproof}

\begin{clm} \label{clm:M_in_H}
There exists a $c_H$-monochromatic component $M'$ such that $V(M) \cap V(G_0) \subseteq V(M')$.
\end{clm}

\begin{claimproof}
Let $a,b$ be vertices in $V(M) \cap V(G_0)$.
To prove this claim, it suffices to show that there exists a $c_H$-monochromatic walk in $H$ between $a$ and $b$.
Since $M$ is connected, there exists a path $Q \subseteq M$ between $a$ and $b$.
Then the walk in $H$ obtained from $Q$ by, for each subpath $P$ of $Q \cap G_e$ between two vertices in $X_e$ (for some $e \in U_E$) internally disjoint from $X_e$, replacing $P$ by the path $\tilde{P}$ mentioned in Claim \ref{clm:path_Gbranch_H}, is a $c_H$-monochromatic walk in $H$ between $a$ and $b$.
\end{claimproof}

\begin{clm} \label{clm:path_from_H_to_G}
For every path $P$ in $H$ with both ends in $V(G_0)$ internally disjoint from $V(G_0)$, there exists a path $\widehat{P}$ in $G$ with $|E(\widehat{P})| \leq (k+2)(\theta+1)^2|E(P)|$ having the same ends as $P$.
\end{clm}

\begin{claimproof}
We may assume that $P$ has an internal vertex, for otherwise we are done by taking $\widehat{P}=P$.
Let $a$ and $b$ be the ends of $P$.
Let $I$ be the set of internal vertices of $P$.
Note that either $I \subseteq S_{e,Y}$ for some $e \in U_E$ and $Y \in \P_e$, or $I \subseteq S_e \cup S'_e$ for some $e \in U_E$. 

We first assume $I \subseteq S_{e,Y}$ for some $e \in U_E$ and $Y \in \P_e$.
Since vertices in $S_{e,Y}$ are pairwise non-adjacent in $H$, $|E(P)|=2$ and $\{a,b\} \subseteq Y$.
By the definition of $\P_e$, there exists a walk in $G$ from $a$ to $b$ with length at most $(2(k+2)(\theta+1)+1)(|X_e|-1) \leq 2(k+2)(\theta+1)^2 = (k+2)(\theta+1)^2|E(P)|$.

So we may assume $I \subseteq S_e \cup S'_e$ for some $e \in U_E$.
By the definition of $H$, there exists a monochromatic component $M_P$ of the unique $L|_{1_L}$-coloring of $G[1_L]$ such that $I \subseteq \{u_{e,M_P},u'_{e,M_P}\}$.
By the definition of $N_H(u_{e,M_P})$ and $N_H(u'_{e,M_P})$, we have $\{a,b\} \subseteq N_{G_e}^{\leq k}[V(M_P) \cap X_{T_e}-X_e]$.
Hence for every $z \in \{a,b\}$, there exists $z' \in V(M_P) \cap X_{T_e}-X_e$ such that $\dist_G(z,z') \leq \dist_{G_e}(z,z') \leq k$.
Since $L$ is $(T,\X,m,s,r,k)$-legitimate, $V(M_P)$ has weak diameter in $G$ at most $k$.
Since $a'$ and $b'$ are contained in $M_P$, $\dist_G(a',b') \leq k$.
Hence $\dist_G(a,b) \leq \dist_G(a,a')+\dist_G(a',b')+\dist_G(b',b) \leq 3k$.
So there exists a path $\widehat{P}$ in $G$ between $a$ and $b$ with $|E(\widehat{P})| \leq 3k \leq (k+2)(\theta+1)^2|E(P)|$ (since $\theta \geq 1$ by Claim \ref{cl:1}).
\end{claimproof}

\begin{clm} \label{clm:central_dist}
For any $a,b \in V(M) \cap V(G_0)$, $$\dist_G(a,b) \leq (k+2)(\theta+1)^2 \cdot f_3(f_1(f^*(\eta-1,N_1+f_2(N)))).$$
\end{clm}

\begin{claimproof}
By Claim \ref{clm:M_in_H}, $a$ and $b$ are contained in the same $c_H$-monochromatic component $M'$ in $H$.
By Claim \ref{claim4}, $V(M')$ has weak diameter in $H$ at most $f_3(f_1(f^*(\eta-1,N_1+f_2(N))))$.
So there exists a path $Q$ in $H$ between $a$ and $b$ with $|E(Q)| \leq f_3(f_1(f^*(\eta-1,N_1+f_2(N))))$.
By replacing each maximal subpath $P$ of $Q$ between $V(G_0)$ internally disjoint from $V(G_0)$ with the path $\widehat{P}$ mentioned in Claim \ref{clm:path_from_H_to_G}, there exists a walk in $G$ from $a$ to $b$ with length at most $(k+2)(\theta+1)^2 \cdot |E(Q)|$.
Therefore, $\dist_G(a,b) \leq (k+2)(\theta+1)^2 \cdot f_3(f_1(f^*(\eta-1,N_1+f_2(N))))$.
\end{claimproof}

\begin{clm} \label{clm:go_up}
For every $v \in V(M)$, there exists $\iota(v) \in V(M) \cap V(G_0)$ such that $\dist_G(v,\iota(v)) \leq 2(k+2)(\theta+1)^2$.
\end{clm}

\begin{claimproof}
We may assume $v \in V(M)-V(G_0)$, for otherwise we are done by taking $\iota(v)=v$.
So there exists $e \in B$ such that $v \in X_{T_e}-X_e$.
By Claim \ref{claim_G0_and_others}, $V(M) \cap V(G_0) \neq \emptyset$.
So there exists a path $P_v$ in $M \cap G_e$ from $v$ to $X_e$ internally disjoint from $X_e$.

Let $a$ be the end of $P_v$ in $X_e$.
Let $b$ be the neighbor of $a$ in $P_v$.
So $a \in Y_b$.
By Claim \ref{clm:path_same_Y}, $a \in Y_b = Y_v$.

By Claim \ref{claim_G0_and_others}, $v \in Z_{\theta+1} \cap X_{T_e}-X_e$.
So there exists a path $Q$ in $G_e$ from $v$ to $X_e$ internally disjoint from $X_e$ with length at most $(k+2)(\theta+1)$.
Let $q$ be the end of $Q$ in $X_e$.
Note that $q \in Y_v$ by the definition of $Y_v$.

Since both $a$ and $q$ are contained in $Y_v$, by the definition of $\P_e$, there exists a path in $G_e$ between $a$ and $q$ with length at most $(2(k+2)(\theta+1)+1)(|X_e|-1) \leq (2(k+2)(\theta+1)+1)(\theta-1)$.
So $\dist_G(v,a) \leq \dist_G(v,q)+\dist_G(q,a) \leq (k+2)(\theta+1) + (2(k+2)(\theta+1)+1)(\theta-1) \leq 2(k+2)(\theta+1)^2$.
Since $a \in V(P_v) \cap X_e \subseteq V(M) \cap V(G_0)$, we are done by choosing $\iota(v)=a$.
\end{claimproof}

Since the weak diameter in $G$ of $V(M)$ is greater than $f^*(\eta,N)$, there exist vertices $a$ and $b$ in $V(M)$ such that $\dist_G(a,b)>f^*(\eta,N)$.
By Claims \ref{clm:central_dist} and \ref{clm:go_up}, 
	\begin{align*}
		& \dist_G(a,b) \\
		\leq & \dist_G(a,\iota(a)) + \dist_G(\iota(a),\iota(b)) + \dist_G(\iota(b),b) \\
		\leq & 2(k+2)(\theta+1)^2 + (k+2)(\theta+1)^2 \cdot f_3(f_1(f^*(\eta-1,N_1+f_2(N)))) + 2(k+2)(\theta+1)^2 \\
		= & (k+2)(\theta+1)^2 \cdot (4+f_3(f_1(f^*(\eta-1,N_1+f_2(N))))) \\
		\leq & f^*(\eta,N),
	\end{align*}
a contradiction.
This proves the lemma.
\end{pf}

\section{Coloring graphs with an excluded minor} \label{sec:main_color}

Now we are ready to prove Theorem \ref{tw_intro}, which follows from the following theorem by the Grid Minor Theorem \cite[(1.5)]{rs_V}.

\begin{theorem} \label{tw_coloring}
For any positive integers $w$ and $k$, there exists a positive integer $N$ such that for every graph $G$ of tree-width at most $w$, for every 2-list-assignment $L$ of $G$, and for every pre-$L$-coloring $c_0$ with weak diameter in $G$ at most $k$, $c_0$ can be extended to an $L$-coloring of $G$ with weak diameter in $G$ at most $N$. 
\end{theorem}

\begin{pf}
Let $w$ and $k$ be positive integers.
Let $N_1$ be the integer $N^*$ mentioned in Lemma \ref{all_centered} by taking $(k,r)=(w+1,2)$.
Let $N_2=\max\{N_1,4\}$.
Let $f$ be the function $f^*$ mentioned in Lemma \ref{tree_extension_list} by taking $(\theta,s,r,k)=(w+1,w+1,1,k)$. 
Define $N=f(w+1,N_2)$.

Let $G,L,c_0$ be as defined in the theorem.

Since $G$ has tree-width at most $w$, there exists a rooted tree-decomposition $(T,\X=(X_t: t \in V(T)))$ of $G$ with width at most $w$.
We may assume that the bag of the root of $T$ is non-empty.
So $(T,\X)$ is a $(w+1,w+1)$-construction of $G$.

Let $S$ be the domain of $c_0$.
For every $v \in S$, let $L'(v)=\{c_0(v)\}$.
For every $v \in V(G)-S$, let $L'(v)=L(v)$.
So $L'$ is a list-assignment of $G$ such that $1_{L'}=S$ and $|L'(v)| \in \{1,2\}$ for every $v \in V(G)$.
For every $t \in V(T)$, $1_{L'} \cap X_t$ is $(|X_t|,0)$-centered in $G[X_t]$ and hence $(w+1,1)$-centered in $G[X_t]$.
Since $1_{L'}=S$, the unique $L'_{1_{L'}}$-coloring of $G[1_{L'}]$ is $c_0$ and has weak diameter in $G$ at most $k$.

Hence $L'$ is a $(T,\X,2,w+1,1,k)$-legitimate list-assignment of $G$.
We shall show that $L'$ is hereditarily $(T,\X,2,w+1,1,k,N_2)$-extendable.

Let $H$ be an induced subgraph of $G$.
Let $L_H$ be a $(T,\X \cap V(H),2,w+1,1,k)$-legitimate list-assignment of $H$ such that $L_H(v) \supseteq L'(v)$ for every $v \in V(H)$.
Let $t \in V(T)$.
Let $G_t$ be a child-extension at $t$ with respect to $(T,\X \cap V(H))$.
Let $L_t$ be a list-assignment of $G_t$ with $L_t|_{X_t \cap V(H)}=L'|_{X_t \cap V(H)}$.
Note that every component of $G_t$ disjoint from $X_t \cap V(H)$ contains at most two vertices; every component of $G_t$ intersecting $X_t \cap V(H)$ has its vertex-set in $N_{G_t}^{\leq 2}[X_t \cap V(G_t)]$.
So the vertex-set of each component of $G_t$ is $(w+1,2)$-centered in $G_t$.
Hence by Lemma \ref{all_centered}, there exists a $L_t$-coloring of $G_t$ with weak diameter in $G_t$ at most $N_1 \leq N_2$.

This shows that $L'$ is hereditarily $(T,\X,2,w+1,1,k,N_2)$-extendable.
By Lemma \ref{tree_extension_list}, every coloring of the empty set can be extended to an $L'$-coloring $c$ of $G$ with weak diameter in $G$ at most $N$.
Note that $c$ is an $L$-coloring that can be obtained by extending $c_0$.
This proves the theorem.
\end{pf}

\bigskip

Our next goal in this section is to prove Theorem \ref{minor_intro}.
We first point out a folklore result about tree-width of graphs with bounded Euler genus and bounded radius (Corollary \ref{radius_tw}).

We first recall notions related to layered tree-width.
A {\it layering} of a graph $G$ is an ordered partition $(V_1,V_2,...)$ of $V(G)$ into possibly empty sets such that for every edge $e$ of $G$, there exists $i \in {\mathbb N}$ such that the ends of $e$ are contained in $V_i \cup V_{i+1}$.
The {\it layered tree-width} of a graph $G$ is the minimum integer $w$ such that there exists a tree-decomposition $(T,(X_t: t \in V(T)))$ of $G$ and a layering $(V_1,V_2,...)$ of $G$ such that $|X_t \cap V_i| \leq w$ for every $t \in V(T)$ and $i \in {\mathbb N}$.

\begin{theorem}[{{\cite[Theorem 12]{dmw}}}] \label{genus_ltw}
For every integer $g$, every graph that can be drawn in a surface of Euler genus $g$ has layered tree-width at most $2g+3$.
\end{theorem}

\begin{corollary} \label{radius_tw}
For any nonnegative integers $g,d$, if $G$ is a graph such that for every component $C$ of $G$, $C$ can be drawn in a surface of Euler genus at most $g$ and $V(C)$ is $(1,d)$-centered in $G$, then the tree-width of $G$ is at most $(2g+3)(2d+1)-1$.
\end{corollary}

\begin{pf}
We may assume that $G$ is connected since the tree-width of a graph equals the maximum of the tree-width of its components.
So $V(G)$ is $(1,d)$-centered in $G$.
Hence there exists $v \in V(G)$ such that $V(G) \subseteq N_G^{\leq d}[\{v\}]$.
By Theorem \ref{genus_ltw}, there exists a tree-decomposition $(T,(X_t: t \in V(T)))$ of $G$ and a layering $(V_1,V_2,...)$ of $G$ such that $|X_t \cap V_i| \leq 2g+3$ for every $t \in V(T)$ and $i \in {\mathbb N}$.
Since $G$ is connected and $V(G) \subseteq N_G^{\leq d}[\{v\}]$, there are at most $2d+1$ indices $i$ with $V_i \neq \emptyset$.
So $|X_t| \leq (2g+3)(2d+1)$ for every $t \in V(T)$.
\end{pf}

\bigskip

Next, we show that near embeddings are 3-choosable with bounded weak diameter, which is a step toward proving Theorem \ref{minor_intro}.

Let $\Sigma$ be a surface. 
Let $p$ be a nonnegative integer.
We say that a graph $G$ is {\it $p$-nearly embeddable in $\Sigma$} if the following statements hold.
	\begin{itemize}
		\item There exist edge-disjoint subgraphs $G_0,G_1,...,G_p$ of $G$ such that $G=\bigcup_{i=0}^pG_i$.
		\item For each $i \in [p]$, there exists a cyclic ordering $\Omega_i$ of $V(G_0) \cap V(G_i)$.
		\item There exist pairwise disjoint closed disks $\Delta_1,\Delta_2,...,\Delta_p$ in $\Sigma$ such that $G_0$ can be drawn in $\Sigma$ in a way that the drawing only intersects $\bigcup_{i=1}^p\Delta_i$ at vertices of $G_0$, no vertex of $G_0$ is contained in the interior of $\bigcup_{i=1}^p\Delta_i$, and for each $i \in [p]$, the set of vertices of $G_0$ contained in the boundary of $\Delta_i$ is $V(G_0) \cap V(G_i)$, and the ordering of the vertices in $V(G_0) \cap V(G_i)$ appearing in the boundary of $\Delta_i$ equals $\Omega_i$.
		\item For each $i \in [p]$, if we denote $\Omega_i$ by $(v_1,v_2,...,v_{k_i})$ for some integer $k_i$, then there exists a path-decomposition $(P_i,\X_i)$ of $G_i$ such that $\lvert V(P_i) \rvert =k_i$, and for each $j \in [k_i]$, the bag at the $j$-th vertex of $P_i$ contains $v_j$ and has size at most $p$.
	\end{itemize}
We call $(G_0,G_1,...,G_p,\Delta_1,...,\Delta_p,\Omega_1,...,\Omega_p,((P_i,\X_i): i \in [p]))$ a {\it witness} of $G$.

A {\it geodesic} in a graph $G$ is a path $P$ in $G$ such that $|E(P)|=\dist_G(a,b)$, where $a$ and $b$ are the ends of $P$.

\begin{lemma} \label{color_near_embedding_1}
For every surface $\Sigma$ and every nonnegative integer $p$, there exists a positive integer $N$ such that every graph that is $p$-nearly embeddable in $\Sigma$ is 3-choosable with weak diameter in $G$ at most $N$. 
\end{lemma}

\begin{pf}
Let $\Sigma$ be a surface.
Let $p$ be a nonnegative integer.
Let $N_1$ be the integer $N$ mentioned in Theorem \ref{surface_ch} by taking $\Sigma=\Sigma$.
Let $g$ be the Euler genus of $\Sigma$.
Define $N$ to be the integer $N$ mentioned in Theorem \ref{tw_coloring} by taking $(w,k)=((2g+3)(4(2p+1)(N_1+1)+1)p,N_1)$.
Note that $N \geq N_1$ by the statement of Theorem \ref{tw_coloring}.

Let $G$ be a graph that is $p$-nearly embeddable in $\Sigma$.
Let $L$ be a 3-list-assignment of $G$.
It suffices to show that there exists an $L$-coloring of $G$ with weak diameter in $G$ at most $N$.

Let $(G_0,G_1,...,G_p,\Delta_1,...,\Delta_p,\Omega_1,...,\Omega_p,((P_i,\X_i): i \in [p]))$ be a witness of $G$.
Let $G_0' = G_0-\bigcup_{i=1}^pV(G_i)$.
Note that $G_0'=G-\bigcup_{i=1}^pV(G_i)$, so $G_0'$ is an induced subgraph of $G$, and $G_0'$ can be drawn in $\Sigma$.
By Theorem \ref{surface_ch}, there exists an $L|_{V(G_0')}$-coloring $c_0$ of $G_0'$ with weak diameter in $G_0'$ at most $N_1$.
For every $i \in [p]$, 
	\begin{itemize}
		\item let $D_i = V(G_0') \cap N_G^{\leq 1}[V(G_0) \cap V(G_i)]$,
		\item let $\F_i=\{M: M$ is a $c_0$-monochromatic component in $G_0'$ intersecting $D_i\}$, and
		\item for every $M \in \F_i$ and $\{a,b\} \subseteq V(M)$, let $P_{ab}$ be a geodesic in $G_0'$ between $a$ and $b$.
	\end{itemize}
Note that $|E(P_{ab})| \leq N_1$ since $V(M)$ has weak diameter in $G_0'$ at most $N_1$.
For every $i \in [p]$, let $$W_i = \bigcup_{M \in \F_i}(M \cup \bigcup_{a,b \in V(M)}P_{ab}).$$
Note that $V(W_i) = \bigcup_{M \in \F_i}(V(M) \cup \bigcup_{a,b \in V(M)}V(P_{ab})) \subseteq N_{G_0'[V(W_i)]}^{\leq N_1}[\bigcup_{M \in \F_i}V(M)]$.
For every $i \in [p]$, since each $M \in \F_i$ intersects $D_i$ and the geodesics realizing the weak diameter of $V(M)$ are included in $W_i$, we have $\bigcup_{M \in \F_i}V(M) \subseteq N_{G_0'[V(W_i)]}^{\leq N_1}[D_i]$, so $V(W_i) \subseteq N_{G_0'[V(W_i)]}^{\leq 2N_1}[D_i]$.
Since $G_0'$ is a subgraph of $G_0$, $$V(W_i) \subseteq N_{G_0[V(W_i)]}^{\leq 2N_1}[D_i]$$ for every $i \in [p]$.

Let $H_0$ be the graph obtained from $G_0$ by, for each $i \in [p]$, adding a new vertex $v_i$ drawn in the interior of $\Delta_i$ adjacent to all vertices in $V(G_0) \cap V(G_i)$ and adding edges such that there exists a cycle $C_i$ in $H_0$ with $V(C_i)=V(G_0) \cap V(G_i)$ passing through its vertices in the order same as $\Omega_i$.

\medskip

\noindent{\bf Claim 1:} For every $i \in [p]$, $V(W_i) \cup (V(G_0) \cap V(G_i)) \cup \{v_i\}$ is $(1,2N_1+2)$-centered in $H_0$.

\noindent{\bf Proof of Claim 1:}
Since $H_0$ is a supergraph of $G_0$, $$V(W_i) \subseteq N_{G_0[V(W_i)]}^{\leq 2N_1}[D_i] \subseteq N_{H_0[V(W_i)]}^{\leq 2N_1}[D_i]$$ for every $i \in [p]$. 
For every $i \in [p]$, since $D_i$ is contained in the closed neighborhood of $V(G_0) \cap V(G_i)$, we know $$N_{H_0[V(W_i)]}^{\leq 2N_1}[D_i] \subseteq N_{H_0[V(W_i) \cup (V(G_0) \cap V(G_i))]}^{\leq 2N_1+1}[V(G_0) \cap V(G_i)];$$ since $v_i$ is adjacent to all vertices in $V(G_0) \cap V(G_i)$, we have $$N_{H_0[V(W_i) \cup (V(G_0) \cap V(G_i))]}^{\leq 2N_1+1}[V(G_0) \cap V(G_i)] \subseteq N_{H_0[V(W_i) \cup (V(G_0) \cap V(G_i)) \cup \{v_i\}]}^{\leq 2N_1+2}[\{v_i\}].$$
Therefore, for every $i \in [p]$, $$V(W_i) \subseteq N_{H_0[V(W_i) \cup (V(G_0) \cap V(G_i)) \cup \{v_i\}]}^{\leq 2N_1+2}[\{v_i\}]$$ is $(1,2N_1+2)$-centered in $H_0$, which implies that $V(W_i) \cup (V(G_0) \cap V(G_i)) \cup \{v_i\}$ is $(1,2N_1+2)$-centered in $H_0$.
$\Box$

\medskip

Let $$H_1 = H_0[\bigcup_{i=1}^p(V(W_i) \cup (V(G_0) \cap V(G_i)) \cup \{v_i\})].$$
By Claim 1, $V(W_i) \cup (V(G_0) \cap V(G_i)) \cup \{v_i\}$ has weak diameter in $H_0$ at most $2(2N_1+2)$ for each $i \in [p]$.
It together with Claim 1 imply that for every component of $H_1$, its vertex-set is $(1,(2p+1)(2N_1+2))$-centered in $H_1$.
By Corollary \ref{radius_tw}, $H_1$ has tree-width at most $(2g+3)(2(2p+1)(2N_1+2)+1)-1$.

Let $(T^1,\Y^1)$ be a tree-decomposition of $H_1$ with width at most $(2g+3)(2(2p+1)(2N_1+2)+1)-1$.
Denote $\Y^1$ by $(Y^1_t: t \in V(T^1))$.

By the definition of a witness, for every $i \in [p]$, there exists a bijection $\iota_i$ from $V(G_0) \cap V(G_i)$ to $V(P_i)$ such that for every $v \in V(G_0) \cap V(G_i)$, the bag at $\iota_i(v)$ in $(P_i,\X_i)$ contains $v$.
For every $i \in [p]$, denote the bag at $\iota_i(v)$ in $(P_i,\X_i)$ by $X_{i,v}$.

Let $H_2=H_1 \cup \bigcup_{i=1}^pG_i$.
For every $t \in V(T^1)$, let $$Y^2_t = Y^1_t \cup \bigcup_{i=1}^p\bigcup_{v \in Y^1_t \cap (V(G_0) \cap V(G_i))}X_{i,v}.$$
Since every bag of $(P_i,\X_i)$ has size at most $p$, $|Y^2_t| \leq |Y^1_t|p \leq (2g+3)(4(2p+1)(N_1+1)+1)p$.
Let $\Y^2=(Y_t^2: t \in V(T^1))$.
By the existence of each cycle $C_i$ passing through $V(G_0) \cap V(G_i)$, it is straightforward to show that $(T^1,\Y^2)$ is a tree-decomposition of $H_2$ with width at most $(2g+3)(4(2p+1)(N_1+1)+1)p$.

Let $H = G[V(H_2) \cap V(G)]$.
Since $H$ is a subgraph of $H_2$, the tree-width of $H$ is at most $(2g+3)(4(2p+1)(N_1+1)+1)p$.
For every $v \in \bigcup_{i=1}^p\bigcup_{M \in \F_i}V(M)$, let $c_H(v)=c_0(v)$.
Note that $H[\bigcup_{i=1}^p\bigcup_{M \in \F_i}V(M)]$ is an induced subgraph of $G_0'$ and hence is an induced subgraph of $G$.
So the set of $c_H$-monochromatic components is $\bigcup_{i=1}^p\F_i$.
Hence $c_H$ is a pre-$L|_{V(H)}$-coloring of $H$ with weak diameter in $H$ at most $N_1$ since $H$ includes those paths $P_{ab}$'s.
By Theorem \ref{tw_coloring}, $c_H$ can be extended to an $L|_{V(H)}$-coloring $c_H^*$ of $H$ with weak diameter in $H$ at most $N$.

For every $v \in V(G_0')$, let $c(v)=c_0(v)$.
For every $v \in V(G)-V(G_0') \subseteq \bigcup_{i=1}^pV(G_i) \subseteq V(H)$, let $c(v)=c_H^*(v)$.
Note that $c$ is an $L$-coloring of $G$.
Let $Q$ be a $c$-monochromatic component in $G$.
To prove this lemma, it suffices to show that $V(Q)$ has weak diameter in $G$ at most $N$.

If $V(Q) \cap \bigcup_{i=1}^p(V(G_i) \cup D_i)=\emptyset$, then $Q$ is a connected subgraph of $G_0'$, so $Q$ is contained in a $c_0$-monochromatic in $G_0'$, and hence $V(Q)$ has weak diameter in $G_0'$ (and hence in $G$) at most $N_1 \leq N$.

So we may assume $V(Q) \cap \bigcup_{i=1}^p(V(G_i) \cup D_i) \neq \emptyset$.
If $V(Q) \cap \bigcup_{i=1}^pD_i =\emptyset$, then since $V(Q) \cap \bigcup_{i=1}^p(V(G_i) \cup D_i) \neq \emptyset$, $Q$ is a connected subgraph of $G[\bigcup_{i=1}^pV(G_i)] \subseteq H$, so $Q$ is contained in a $c_H$-monochromatic component of $H$, and hence $V(Q)$ has weak diameter in $H \subseteq G$ at most $N$.

So we may assume $V(Q) \cap \bigcup_{i=1}^pD_i \neq \emptyset$.
Hence $V(Q) \cap \bigcup_{i=1}^p\bigcup_{M \in \F_i}V(M) \neq \emptyset$.
Since $c|_{V(G_0')}=c_0$, $V(Q) \cap V(G_0') \subseteq \bigcup_{i=1}^p\bigcup_{M \in \F_i}V(M)$.
Hence $V(Q) \subseteq V(H) \subseteq V(G)$ and $Q$ is a connected subgraph of $G[V(H_2) \cap V(G)]=H$.
Since $V(Q) \cap V(G_0') \subseteq \bigcup_{i=1}^p\bigcup_{M \in \F_i}V(M)$, $c_0|_{V(Q) \cap V(G_0')} = c_H|_{V(Q) \cap V(G_0')}$.
Hence 
	\begin{align*}
		c|_{V(Q)} = & c|_{V(Q) \cap V(G_0')} \cup c|_{V(Q)-V(G_0')} \\
		= & c_0|_{V(Q) \cap V(G_0')} \cup c_H^*|_{V(Q)-V(G_0')} \\
		= & c_H|_{V(Q) \cap V(G_0')} \cup c_H^*|_{V(Q)-V(G_0')}.
	\end{align*}
Since $c_H^*$ is obtained by extending $c_H$, $c|_{V(Q)} = c_H^*|_{V(Q)}$, so $V(Q)$ is contained in a $c_H^*$-monochromatic component in $H$.
Hence $V(Q)$ has weak diameter in $H \subseteq G$ at most $N$.
This proves the lemma.
\end{pf}

\begin{lemma} \label{color_near_embedding_2}
For every surface $\Sigma$ and every nonnegative integer $p$, there exists a positive integer $N$ such that if $G$ is a graph such that $G-Z$ is $p$-nearly embeddable in $\Sigma$ for some $Z \subseteq V(G)$ with $|Z| \leq p$, then $G$ is 3-choosable with weak diameter in $G$ at most $N$. 
\end{lemma}

\begin{pf}
Let $\Sigma$ be a surface.
Let $p$ be a nonnegative integer.
Let $N_1$ be the positive integer $N$ mentioned in Lemma \ref{color_near_embedding_1} by taking $(\Sigma,p)=(\Sigma,p)$.
Define $N$ to be the integer $N^*$ mentioned in Lemma \ref{add_centered} by taking $(k,r,N)=(p,0,N_1)$.

Let $G$ be a graph.
Let $Z \subseteq V(G)$ with $|Z| \leq p$ such that $G-Z$ is $p$-nearly embeddable in $\Sigma$.
Let $L$ be a 3-list-assignment of $G$.
It suffices to show that there exists an $L$-coloring of $G$ with weak diameter in $G$ at most $N$.

Since $G-Z$ is $p$-nearly embeddable in $\Sigma$, by Lemma \ref{color_near_embedding_1}, there exists an $L|_{V(G)-Z}$-coloring $c_1$ of $G-Z$ with weak diameter in $G-Z$ (and hence in $G$) at most $N_1$.
Let $c_Z$ be an $L|_Z$-coloring of $G[Z]$.
Since $|Z| \leq p$, $Z$ is $(p,0)$-centered in $G$.
By Lemma \ref{add_centered}, $c_1 \cup c_Z$ is an $L$-coloring of $G$ with weak diameter in $G$ at most $N$.
\end{pf}

\begin{lemma} \label{small_extension}
Let $d$ and $N$ be nonnegative integers.
Let $G$ be a graph.
Let $L$ be a list-assignment of $G$.
Let $S \subseteq V(G)$.
Let $G'$ be the graph obtained from $G[S]$ by for each component $C$ of $G-S$, adding edges such that $N_G^{\leq 1}[V(C)] \cap S$ is a clique in $G'$.
If $G'$ has an $L|_S$-coloring with weak diameter in $G'$ at most $N$, and every component of $G-S$ has diameter at most $d$, then $G$ has an $L$-coloring with weak diameter in $G$ at most $(d+2)N+2d+2$.
\end{lemma}

\begin{pf}
Let $c'$ be an $L|_S$-coloring of $G'$ with weak diameter in $G'$ at most $N$.
Define $c$ to be an arbitrary $L$-coloring of $G$ such that $c|_S=c'$.
Let $M$ be a $c$-monochromatic component in $G$.
To prove this lemma, it suffices to show that the weak diameter of $V(M)$ in $G$ is at most $(d+2)N+2d+2$.

If $V(M) \cap S=\emptyset$, then since $M$ is a connected subgraph of $G$, $M$ is contained in a component of $G-S$, which as diameter at most $d$, so the weak diameter in $G$ of $V(M)$ is at most $d$.
So we may assume $V(M) \cap S \neq \emptyset$.

Let $M'=G'[V(M) \cap S]$.
For every path $P$ in $G$ between two vertices in $S$ internally disjoint from $S$, the internal vertices of $P$ are contained in some component $C$ of $G-S$, so there exists an edge of $G'$ between the ends of $P$.
This implies that $M'$ is a connected subgraph of $G'$ since $M$ is a connected subgraph of $G$.
Hence $M'$ is contained in a $c'$-monochromatic component in $G'$.
Therefore, $V(M')$ has weak diameter in $G'$ at most $N$.

Note that for every edge $xy$ of $G'$, either $xy \in E(G)$ or there exists a component $C_{xy}$ of $G-S$ such that $\{x,y\} \subseteq N_G^{\leq 1}[V(C_{xy})] \cap S$; in either case, there exists a path $P_{xy}$ in $G$ between $x$ and $y$ with length at most $d+2$ since every component of $G-S$ has diameter at most $d$.
So for every path $P$ in $G'$ between two vertices $x$ and $y$, there exists a path $P'$ in $G$ between $x$ and $y$ such that $|E(P')| \leq (d+2)|E(P)|$.
For any $a,b \in V(M) \cap S=V(M')$, since $V(M')$ has weak diameter in $G'$ at most $N$, $\dist_G(a,b) \leq (d+2)\dist_{G'}(a,b) \leq (d+2)N$.

For every $v \in V(M)-S$, there exists a component $C_v$ of $G-S$ with $v \in V(C_v)$, and since $V(M) \cap S \neq \emptyset$, there exists a path in $M$ from $v$ to a vertex $\iota(v)$ in $N_G^{\leq 1}[V(C_v)] \cap S$ internally disjoint from $S$, so $\dist_G(v,\iota(v)) \leq d+1$.
For every $v \in V(M) \cap S$, let $\iota(v)=v$.
Note that $\iota(v) \in V(M) \cap S$ and $\dist_G(v,\iota(v)) \leq d+1$ for every $v \in V(M)$.
Hence for any $a,b \in V(M)$, $\dist_G(a,b) \leq \dist_G(a,\iota(a))+\dist_G(\iota(a),\iota(b))+\dist_G(\iota(b),b) \leq (d+1)+(d+2)N+(d+1) = (d+2)N+2d+2$.
This shows that the weak diameter in $G$ of $V(M)$ is at most $(d+2)N+2d+2$.
\end{pf}

\bigskip

For a tree-decomposition $(T,(X_t: t \in V(T)))$ of a graph $G$ and for $t \in V(T)$, the {\it torso} at $t$ is the graph obtained from $G[X_t]$ by for each neighbor $t'$ of $t$ in $T$, adding edges such that $X_t \cap X_{t'}$ is a clique.

\begin{lemma} \label{color_torso}
For any positive integers $p,N$, there exists a positive integer $N^*$ such that the following holds.
Let $m \geq 2$ be an integer, $G$ a graph, and $L$ an $m$-list-assignment of $G$.
If $(T,\X)$ is a tree-decomposition of $G$ with adhesion at most $p$ such that for every $t \in V(T)$ and every subgraph $W$ of the torso at $t$, there exists an $L|_{V(W)}$-coloring of $W$ with weak diameter in $W$ at most $N$, then there exists an $L$-coloring of $G$ with weak diameter in $G$ at most $N^*$.
\end{lemma}

\begin{pf}
Let $p$ and $N$ be positive integers.
Let $f$ be the function $f^*$ mentioned in Lemma \ref{tree_extension_list} by taking $(\theta,s,r,k)=(p,p,1,1)$.
Define $N^* = f(p,3N+4)$.

Let $m \geq 2$ be an integer.
Let $G$ be a graph.
Let $L$ be an $m$-list-assignment of $G$.
Let $(T,\X=(X_t: t \in V(T)))$ be a tree-decomposition of $G$ with adhesion at most $p$ such that for every $t \in V(T)$ and every subgraph $W$ of the torso at $t$, there exists an $L|_{V(W)}$-coloring of $W$ with weak diameter in $W$ at most $N$.

Let $t_0$ be a node of $T$ with $X_{t_0} \neq \emptyset$.
Let $T^*$ be the tree obtained from $T$ by adding a new node $t^*$ adjacent to $t_0$.
We treat $T^*$ as a rooted tree rooted at $t^*$.
Let $v_0 \in X_{t_0}$.
Let $X^*_{t^*}=\{v_0\}$.
For every $t \in V(T)$, let $X^*_t=X_t$.
Let $\X^*=(X^*_t: t \in V(T^*))$.
Then $(T^*,\X^*)$ is a rooted tree-decomposition of $G$.
Since $p \geq 1=|X^*_{t^*}|$, $(T^*,\X^*)$ is a $(p,p)$-construction of $G$.

We shall show that $L$ is hereditarily $(T^*,\X^*,m,p,1,1,3N+4)$-extendable.
Since $L$ is an $m$-list-assignment with $1_L=\emptyset$, we know that $L$ is $(T^*,\X^*,m,p,1,1)$-legitimate.
So $L$ satisfies (H1).

Let $H$ be an induced subgraph of $G$.
Let $L_H$ be a $(T^*,\X^* \cap V(H),m,p,1,1)$-legitimate list-assignment of $H$ with $L_H(v) \supseteq L(v)$ for every $v \in V(H)$.
Note that $|L_H(v)| \geq |L(v)|=m$ for every $v \in V(H)$.
Since $L_H$ satisfies (L1), $L_H=L|_{V(H)}$.
Let $t \in V(T^*)$.
Let $G_t$ be a child-extension at $t$ with respect to $(T^*,\X^* \cap V(H))$.
Let $L_t$ be a list-assignment of $G_t$ with $L_t|_{X^*_t \cap V(H)}=L_H|_{X^*_t \cap V(H)}$.

To show that $L$ satisfies (H2) for being hereditarily $(T^*,\X^*,m,p,1,1,3N+4)$-extendable, it suffices to show that there exists an $L_t$-coloring of $G_t$ with weak diameter in $G_t$ at most $3N+4$.
If $t=t^*$, then $|X^*_{t^*}|=1$, so every component of $G_t$ has diameter at most 4, so there exists an $L_t$-coloring of $G_t$ with weak diameter in $G_t$ at most $4 \leq 3N+4$.
So we may assume $t \in V(T)$.
Hence $X^*_t=X_t$.
So $G_t$ is also a child-extension at $t$ with respect to $(T,\X \cap V(H))$.
Let $G_t'$ be the graph obtained from $G_t[X^*_t \cap V(H)]=G_t[X_t \cap V(H)]$ by for each component $C$ of $G_t-(X^*_t \cap V(H))$, adding edges such that $N_{G_t}^{\leq 1}[V(C)] \cap (X^*_t \cap V(H))$ is a clique.
Since $G_t$ is a child-extension at $t$ with respect to $(T,\X \cap V(H))$, for each component $C$ of $G_t-(X^*_t \cap V(H))$, there exists a child $t_C$ of $t$ in $(T,\X)$ such that $N_{G_t}^{\leq 1}[V(C)] \cap X^*_t \cap V(H) \subseteq X_t \cap X_{t_C}$.
Hence $G_t'$ is a subgraph of the torso at $t$ with respect to $(T,\X)$.
By the property of $(T,\X)$, there exists an $L|_{X^*_t \cap V(H)}$-coloring of $G_t'$ with weak diameter in $G_t'$ at most $N$.
Since $L_t|_{V(G_t')}=L_t|_{X^*_t \cap V(H)}=L_H|_{X^*_t \cap V(H)}=L|_{X^*_t \cap V(H)}$, there exists an $L_t|_{V(G_t')}$-coloring of $G_t'$ with weak diameter in $G_t'$ at most $N$.
Since $G_t$ is a child-extension at $t$, every component of $G_t-(X^*_t \cap V(H))$ has at most two vertices and hence has diameter at most $1$.
So by Lemma \ref{small_extension}, there exists an $L_t$-coloring of $G_t$ with weak diameter in $G_t$ at most $3N+4$.

Therefore, $L$ is hereditarily $(T^*,\X^*,m,p,1,1,3N+4)$-extendable.
By Lemma \ref{tree_extension_list}, the precoloring on the empty graph can be extended to an $L$-coloring of $G$ with weak diameter in $G$ at most $f(p,3N+4)=N^*$.
This proves the lemma.
\end{pf}

\bigskip

We will use the following structure theorem proved in \cite{rs_XVI}.

\begin{theorem}[{\cite[Theorem 1.3]{rs_XVI}}] \label{rs_structure}
For every graph $H$, there exists a positive integer $p$ such that if $G$ is an $H$-minor free graph, then there exists a tree-decomposition $(T,(X_t: t \in V(T)))$ of $G$ of adhesion at most $p$ such that for every $t \in V(T)$, there exists $Z_t \subseteq X_t$ with $|Z_t| \leq p$ such that the graph obtained from the torso at $t$ by deleting $Z_t$ is $p$-nearly embeddable in a surface in which $H$ cannot be drawn.
\end{theorem}

Now we are ready to prove Theorem \ref{minor_intro}.
The following is a restatement.

\begin{theorem} \label{minor_color}
For every graph $H$, there exists a positive integer $N$ such that for every $H$-minor free graph $G$ and every 3-list-assignment $L$ of $G$, there exists an $L$-coloring of $G$ with weak diameter in $G$ at most $N$.
\end{theorem}

\begin{pf}
Let $H$ be a graph.
Let $p$ be the positive integer $p$ mentioned in Theorem \ref{rs_structure} by taking $H=H$.
For every surface $\Sigma$, let $N_\Sigma$ be the positive integer $N$ mentioned in Lemma \ref{color_near_embedding_2} by taking $(\Sigma,p)=(\Sigma,p)$.
Let $N_1 = \max_\Sigma N_\Sigma$, where the maximum is over all surfaces $\Sigma$ in which $H$ cannot be drawn.
Define $N$ to be the integer $N^*$ mentioned in Lemma \ref{color_torso} by taking $(p,N)=(p,N_1)$.

Let $G$ be an $H$-minor free graph.
Let $L$ be a 3-list-assignment of $G$.

By Theorem \ref{rs_structure}, there exists a tree-decomposition $(T,(X_t: t \in V(T)))$ of $G$ of adhesion at most $p$ such that for every $t \in V(T)$, there exists $Z_t \subseteq X_t$ with $|Z_t| \leq p$ such that the graph obtained from the torso at $t$ by deleting $Z_t$ is $p$-nearly embeddable in a surface $\Sigma_t$ in which $H$ cannot be drawn. 
So for every $t \in V(T)$ and subgraph $W$ of the torso at $t$, $W-(Z_t \cap V(W))$ is $p$-nearly embeddable in $\Sigma_t$, and hence by Lemma \ref{color_near_embedding_2}, there exists an $L|_{V(W)}$-coloring of $W$ with weak diameter in $W$ at most $N_{\Sigma_t} \leq N_1$.
By Lemma \ref{color_torso}, there exists an $L$-coloring of $G$ with weak diameter in $G$ at most $N$.
\end{pf}

\section{Bipartite graphs and odd minors} \label{sec:odd_minors}

We shall prove Theorems \ref{bipartite_intro} and \ref{odd_minor_intro} in this section.

We need the following lemmas to prove Theorem \ref{bipartite_intro}.

\begin{lemma} \label{girth_far}
Let $g \geq 4$ be a positive integer.
Let $G$ be a graph with girth at least $g$.
If $C$ is a cycle in $G$, then there exist vertices $x,y$ in $C$ such that $\dist_G(x,y) \geq \lfloor g/4 \rfloor$.
\end{lemma}

\begin{pf}
Since $g \geq 4$, there exist disjoint subpaths $P_1$ and $P_2$ of $C$ with length $\lfloor g/4 \rfloor$.
Let $P$ be a path in $G$ from $V(P_1)$ to $V(P_2)$ such that $|E(P)|$ is as small as possible. 
Let $a_1,a_2$ be the ends of $P$ in $V(P_1)$ and $V(P_2)$, respectively.
By the minimality of $|E(P)|$, $P$ is internally disjoint from $V(P_1) \cup V(P_2)$ and $|E(P)|=\dist_G(a_1,a_2)$.
We may assume $\dist_G(a_1,a_2) \leq \lfloor g/4 \rfloor-1$, for otherwise we are done by choosing $x$ and $y$ to be $a_1$ and $a_2$.
Hence $|E(P)| \leq \lfloor g/4 \rfloor-1$.

For $i \in [2]$, let $x_i$ and $y_i$ be the ends of $P_i$ such that $\dist_{P_i}(x_i,a_i) \geq \dist_{P_i}(y_i,a_i)$, and let $Q_i$ be the subpath of $P_i$ between $x_i$ and $a_i$.
So $|E(Q_i)|=\dist_{P_i}(x_i,a_i) \geq \frac{1}{2}|E(P_i)|= \frac{1}{2} \lfloor \frac{g}{4} \rfloor$.
Let $Q = Q_1 \cup P \cup Q_2$.
So $|E(Q)| \geq |E(Q_1)|+ |E(Q_2)| \geq \lfloor g/4 \rfloor$.

We may assume $\dist_G(x_1,x_2) \leq \lfloor g/4 \rfloor-1$, for otherwise we are done by choosing $x$ and $y$ to be $x_1$ and $x_2$.
Let $R$ be a path in $G$ between $x_1$ and $x_2$ such that $|E(R)|=\dist_G(x_1,x_2)$.
Note that $|E(R)| = \dist_G(x_1,x_2) \leq \lfloor g/4 \rfloor-1 < |E(Q)|$.
So $R \neq Q$.
Let $H$ be the graph with $V(H)=V(R) \cup V(Q)$ such that $E(H)$ equals the symmetric difference of $E(R)$ and $E(Q)$.
Clearly, $H$ is an Eulerian subgraph of $G$ with $E(H) \neq \emptyset$.
So there exists a cycle $C'$ in $H \subseteq G$.
But $|E(C')| \leq |E(Q)|+|E(R)| \leq |E(Q_1)|+|E(P)|+|E(Q_2)|+|E(R)| \leq |E(P_1)|+(\lfloor g/4 \rfloor-1)+|E(P_2)|+(\lfloor g/4 \rfloor-1) = 4 \cdot \lfloor g/4 \rfloor- 2<g$, contradicting that $G$ has girth at least $g$.
\end{pf}

\begin{lemma} \label{girth_chro_non_list}
Let $d,g,p,k$ be positive integers with $p \geq 3k^{k+1}$ and $g \geq 4$.
If there exists a $d$-regular graph with girth at least $g$ and chromatic number at least $p$, then there exist a bipartite graph $G$ with maximum degree at most $dk^k$ and a $k$-list-assignment $L$ of $G$ such that every $L$-coloring of $G$ has a monochromatic component with weak diameter in $G$ greater than $2 \cdot \lfloor \frac{g}{4} \rfloor-3$. 
\end{lemma}

\begin{pf}
Let $H$ be a $d$-regular graph with girth at least $g$ and chromatic number at least $p$.
Let $Q$ be the graph with $V(Q)=\{q_{v,i}: v \in V(H), i \in [k]\} \cup E(H)$ and $E(Q) = \{q_{v,i}e: v \in V(H), i \in [k], e \in E(H)$ is incident with $v$ in $H\}$.
Note that $Q$ is bipartite graph with a bipartition $\{\{q_{v,i}: v \in V(H), i \in [k]\},E(H)\}$.
For any $v \in V(H)$ and $i \in [k]$, $\deg_Q(q_{v,i}) = \deg_H(v)=d$; for every $e \in E(H)$, $\deg_Q(e) = 2k$.
So the maximum degree of $Q$ is at most $\max\{d,2k\}=d$ since $d \geq \chi(H)-1 \geq p-1 \geq 3k^{k+1}-1 \geq 2k$.

A {\it type} is a sequence $(t_1,t_2,...,t_k)$ such that for every $i \in [k]$, $t_i \in [ik]-[(i-1)k]$. 
Let $T$ be the set of types.
Note that there are $k^k$ different types, so $|T|=k^k$.

For each type $t$, create a copy $Q^t$ of $Q$.
For each vertex $z$ in $Q$, let $z^t$ be the copy of $z$ in $Q^t$.
Define $G$ to be the graph obtained from the disjoint union $\bigcup_{t \in T}Q^t$ of all copies by, for each $v \in V(H)$ and $i \in [k]$, identifying all vertices in $\{q_{v,i}^t: t \in T\}$ into a single vertex $q_{v,i}^*$.

Note that $G$ is a bipartite graph, where one set in the bipartition is $\{q^*_{v,i}: v \in V(H), i \in [k]\}$ and the other set is $\{e^t: e \in E(H), t \in T\}$.
Since the maximum degree of $Q$ is at most $d$ and there are $k^k$ different types, the maximum degree of $G$ is at most $dk^k$.

Define a list-assignment $L$ of $G$ as follows:
    \begin{itemize}
        \item for any $v \in V(H)$ and $i \in [k]$, define $L(q^*_{v,i}) = [ik]-[(i-1)k]$, and
        \item for any other vertex $z$ of $G$, $z=e^t$ for some $e \in E(H)$ and type $t=(t_1,t_2,...,t_k)$, and we define $L(z)=\{t_j: j \in [k]\}$.
    \end{itemize}
Note that $L$ is a $k$-list-assignment of $G$.

Let $c$ be an $L$-coloring of $G$.
To prove this lemma, it suffices to show that some $c$-monochromatic component in $G$ has weak diameter in $G$ greater than $2 \cdot \lfloor \frac{g}{4} \rfloor-3$.

For every $v \in V(H)$, let $c_H(v)=(c(q^*_{v,1}),c(q^*_{v,2}),...,c(q^*_{v,k}))$.
Note that $c_H$ is a function from $V(H)$ to $T$.
Since $p \leq \chi(H) \leq \sum_{t \in T}\chi(H[c_H^{-1}(t)])$, there exists $t^* \in T$ such that $\chi(H[c_H^{-1}(t^*)]) \geq \frac{p}{|T|} \geq \frac{3k^{k+1}}{k^k}=3k$.
So there exists a subgraph $H_1$ of $H[c_H^{-1}(t^*)]$ with minimum degree at least $\chi(H[c_H^{-1}(t^*)])-1 \geq 3k-1 \geq 2k$.
Hence $|E(H_1)| \geq k|V(H_1)|$.

Note that the sets $L(e^{t^*})$ (for all $e \in E(H_1)$) are identical sets with size $k$.
So there exists $X \subseteq E(H_1)$ with $|X| \geq |E(H_1)|/k \geq |V(H_1)|$ such that $c(x^{t^*})$ are identical for all $x \in X$.
Hence there exists a cycle $C$ in $H_1 \subseteq H$ with $E(C) \subseteq X$.

Note that for every $x \in X$, $c(x^{t^*})$ is an entry of $t^*$ by the definition of $L(x^{t^*})$, so for every $v \in V(H_1)$, there exists $i \in [k]$ such that $c(q^*_{v,i})=c(x^{t^*})$ by the definition of $c_H$ and $H_1$.
Together with the fact that $C$ is a cycle with $E(C) \subseteq X$, we know that there exists a $c$-monochromatic component $M$ containing $\{x^{t^*}: x \in X\}$.

Suppose to the contrary that $V(M)$ has weak diameter in $G$ at most $2 \cdot \lfloor \frac{g}{4} \rfloor-3$.
Then for any $x,y \in X$, there exists a path $P_{xy}$ in $G$ from $x^{t^*}$ to $y^{t^*}$ with length at most $2 \cdot \lfloor \frac{g}{4} \rfloor-3$; note that we may assume $P_{xy} \subseteq Q^{t^*}$ since each $Q^t$ is a copy of $Q$; and for any end $z_x$ of $x$ and any end $z_y$ of $y$, there exists a walk $W_{z_xz_y}$ in $H$ from $z_x$ to $z_y$ such that the edges in $W_{z_xz_y}$ are the vertices in $V(P_{xy}) \cap \{e^{t^*}: e \in E(H)\}$, so the number of edges in $W_{z_xz_y}$ is at most $\lceil \frac{1}{2}|V(P_{xy})| \rceil \leq \lfloor \frac{g}{4} \rfloor-1$.
Hence for any two vertices $u,v$ in $C$, $\dist_H(u,v) \leq \lfloor \frac{g}{4} \rfloor-1$ since $E(C) \subseteq X$.
But $C \subseteq H$ and $H$ has girth at least $g$, contradicting Lemma \ref{girth_far}.
\end{pf}

\begin{theorem}[{\cite[Theorem 5.13]{m}}] \label{expanders}
Let $q$ be a power of 2.
Then there exist infinitely many positive integers $t$ such that there exists a $(q+1)$-regular graph $G_{q,t}$ on $q^{3t}-q^t$ vertices with girth at least $\log_q|V(G_{q,t})|$ and with chromatic number at least $\frac{q+1}{2\sqrt{q}}+1$.
\end{theorem}

Theorem \ref{bipartite_intro} follows from the following theorem by taking $\Delta = q^{1.5}$.

\begin{theorem}
For every integer $q$ that is a power of 2 with $q \geq 1296$, there exist infinitely many positive integers $t$ such that some graph $G$ with maximum degree at most $q^{1.5}$ is not $\lfloor \frac{\log q}{8\log\log q} \rfloor$-choosable with weak diameter in $G$ at most $t$.
\end{theorem}

\begin{pf}
Let $q$ be a power of 2 with $q \geq 1296$.
Let $t$ be a positive integer such that $t \geq 2$ and the graph $G_{q,t}$ mentioned in Theorem \ref{expanders} exists.
Note that there are infinitely many such integers $t$ by Theorem \ref{expanders}.

Let $k = \lfloor \frac{\log q}{8\log\log q} \rfloor$.
Note that $2k\log k \leq 2 \cdot \frac{\log q}{8\log\log q} \cdot (\log\log q - \log 8 - \log\log\log q) \leq \frac{\log q}{4}$.
Since $q \geq 1296$, $\frac{\log q}{4} \leq \frac{\log q}{2} - \log 6$.
So $k^{k+1} \leq k^{2k} \leq \sqrt{q}/6$. 

Let $p=\frac{q+1}{2\sqrt{q}}+1$.
Since $k^{k+1} \leq \sqrt{q}/6$, $p \geq \sqrt{q}/2 \geq 3k^{k+1}$.
By Theorem \ref{expanders}, there exists a $(q+1)$-regular graph $G_{q,t}$ on $q^{3t}-q^t$ vertices with girth at least $\log_q|V(G_{q,t})|$ and with chromatic number at least $p$.
Since $q \geq 4$, $q^{3t}-q^t \geq q^{2.5t}$, so the girth of $G_{q,t}$ is at least $\log_q|V(G_{q,t})| \geq \log_q q^{2.5t} = 2.5t \geq 5$.
Since $p \geq 3k^{k+1}$, by Lemma \ref{girth_chro_non_list}, there exists a graph $G$ with maximum degree at most $(q+1)k^k \leq 2qk^{k+1} \leq q^{1.5}$ and a $k$-list-assignment $L$ of $G$ such that every $L$-coloring of $G$ has a monochromatic component with weak diameter in $G$ greater than $2 \cdot \lfloor \frac{2.5t}{4} \rfloor-3 \geq t$ for all sufficiently large $t$.
Hence $G$ is not is not $k$-choosable with weak diameter in $G$ at most $t$.
\end{pf}

\bigskip

Finally, we prove Theorem \ref{odd_minor_intro}.
We will use the following folklore structure theorem that can be obtained by combining results in \cite{ggrsv,rs_XVI} and is explicitly stated in \cite{l_conflict}.

\begin{theorem}[{\cite[Theorem 4.4]{l_conflict}}] \label{structure_odd_minor}
For every graph $H$, there exist positive integers $r$ and $\xi$ such that every odd $H$-minor free graph $G$ has a tree-decomposition $(T,(X_t: t \in V(T)))$ such that for every $t \in V(T)$, if $R_t$ is the torso at $t$, then either $R_t$ is $K_r$-minor free, or there exists $Z_t \subseteq X_t$ with $|Z_t| \leq \xi$ such that $R_t-Z_t$ is bipartite.
\end{theorem}

The following is a restatement of Theorem \ref{odd_minor_intro}.

\begin{theorem}
For every graph $H$, there exists a positive integer $N$ such that every odd $H$-minor free graph has a 3-coloring with weak diameter in $G$ at most $N$.
\end{theorem}

\begin{pf}
Let $H$ be a graph.
Let $r$ and $\xi$ be the positive integers $r$ and $\xi$, respectively, mentioned in Theorem \ref{structure_odd_minor} by taking $H=H$.
Let $\theta = r + \xi +2$.
Let $N_1$ be the positive integer $N$ mentioned in Theorem \ref{minor_color} by taking $H=K_r$.
Let $N_2$ be the positive integer $N$ mentioned in Lemma \ref{add_centered} by taking $(k,r,N)=(\xi,0,1)$.
Define $N$ to be the positive integer $N^*$ mentioned in Lemma \ref{color_torso} by taking $(p,N)=(\theta,N_1+N_2)$.

Let $G$ be an odd $H$-minor free graph.
Let $L$ be the list-assignment of $G$ such that $L(v)=[3]$ for every $v \in V(G)$.
It suffices to show that $G$ has an $L$-coloring with weak diameter in $G$ at most $N$.

By Theorem \ref{structure_odd_minor}, there exists a tree-decomposition $(T,(X_t: t \in V(T)))$ of $G$ such that for every $t \in V(T)$, if $R_t$ is the torso at $t$, then either $R_t$ is $K_r$-minor free, or there exists $Z_t \subseteq X_t$ with $|Z_t| \leq \xi$ such that $R_t-Z_t$ is bipartite.
If the adhesion of $(T,(X_t: t \in V(T)))$ is greater than $\theta$, then some torso contains a clique with size at least $\theta+1 \geq r+\xi+3$, so this torso cannot be $K_r$-minor free and cannot be made bipartite by deleting at most $\xi$ vertices.
So the adhesion of $(T,(X_t: t \in V(T)))$ is at most $\theta$.

Let $t \in V(T)$.
Let $R_t$ be the torso at $t$.
Let $W$ be a subgraph of $R_t$.
By Lemma \ref{color_torso}, to prove this theorem, it suffices to show that there exists an $L|_{V(W)}$-coloring of $W$ with weak diameter in $W$ at most $N_1+N_2$.

We first assume that $R_t$ is $K_r$-minor free.
Since $W$ is a subgraph of $R_t$, $W$ is $K_r$-minor free.
By Theorem \ref{minor_color}, $W$ has an $L|_{V(W)}$-coloring with weak diameter in $W$ at most $N_1$.

So we may assume that $R_t$ is not $K_r$-minor free.
Hence by the property of $(T,(X_t: t \in V(T)))$, there exists $Z_t \subseteq X_t$ with $|Z_t| \leq \xi$ such that $R_t-Z_t$ is bipartite.
Since $W$ is a subgraph of $R_t$, $W-Z_t$ is bipartite.
Since $L(v)=[3]$ for every $v \in V(G)$, there exists a proper $L|_{V(W)-Z_t}$-coloring $c_1$ of $W-Z_t$.
So $c_1$ is an $L|_{V(W)-Z_t}$-coloring of $W-Z_t$ with weak diameter in $W$ at most $0$.
Since $|Z_t| \leq \xi$, $Z_t \cap V(W)$ is $(\xi,0)$-centered in $W$.
Let $c_Z$ be a function from $Z_t \cap V(W)$ to $[3]$.
By Lemma \ref{add_centered}, $c_0 \cup c_Z$ is an $L|_{V(W)}$-coloring of $W$ with weak diameter in $W$ at most $N_2$.
This proves the theorem.
\end{pf}

\section{Concluding remarks} \label{sec:concluding}

In this paper we prove that for every graph $H$, every $H$-minor free graph is 3-choosable with bounded weak diameter (Theorem \ref{minor_intro}).
This result strengthens a result of Dvo\v{r}\'{a}k and Norin \cite{dn} for graphs with bounded Euler genus.
Note that the result of Dvo\v{r}\'{a}k and Norin implies that every graph with bounded Euler genus and bounded maximum degree is 3-choosable with bounded clustering. 
The non-list-coloring version was proved by Esperet and Joret \cite{ej}, and several new proofs of this result were obtained recently by proving stronger results \cite{bbeglps,bbegps,demww,l_asdim,lw}.
Some of those proofs show that graphs with bounded layered tree-width and with bounded maximum degree are 3-colorable with bounded clustering and hence can be applied to graphs with certain geometric properties (see \cite{dmw,lw}).
However, none of those proofs works for list-coloring. 

\begin{conjecture} \label{conj_choosable_layer_deg}
For any positive integers $w$ and $\Delta$, there exists a positive integer $N$ such that every graph with layered tree-width at most $w$ and with maximum degree at most $\Delta$ is 3-choosable with clustering at most $N$.
\end{conjecture}

\begin{conjecture}
For every positive integer $w$, there exists a positive integer $N$ such that every graph with layered tree-width at most $w$ is 3-choosable with weak diameter in $G$ at most $N$.
\end{conjecture}

In this paper, we also prove that bipartite graphs with maximum degree at most $\Delta$ is $k$-choosable with bounded weak diameter only when $k=\Omega(\frac{\log\Delta}{\log\log\Delta})$ (Theorem \ref{bipartite_intro}).
It is equivalent to the result that bipartite graphs with maximum degree at most $\Delta$ is $k$-choosable with bounded clustering only when $k=\Omega(\frac{\log\Delta}{\log\log\Delta})$.

\begin{question}
Determine the slowest growing function $f$ such that for every positive integer $\Delta$, there exists a positive integer $N_\Delta$ such that every bipartite graph with maximum degree at most $\Delta$ is $f(\Delta)$-choosable with clustering at most $N_\Delta$. 
\end{question}

Molloy \cite{m_color} proved that triangle-free graphs with maximum degree at most $\Delta$ are $(1+o(1))\frac{\Delta}{\ln\Delta}$-choosble.
So the function $f$ in the above question is $O(\frac{\Delta}{\log\Delta})$ and $\Omega(\frac{\log\Delta}{\log\log\Delta})$.
Recall that Alon and Krivelevich \cite{ak} conjectured that bipartite graphs with maximum degree at most $\Delta$ are $O(\log\Delta)$-choosable; if this conjecture is true, then the function $f$ in the above question is $O(\log\Delta)$.

\bigskip

\bigskip

\noindent{\bf Acknowledgement:}
The second author thanks Zden\v{e}k Dvo\v{r}\'{a}k and Sergey Norin for bringing \cite{dn} to his attention and for discussions that motivated this work.

\end{document}